\documentclass[12pt]{amsart}
\usepackage{amssymb}
\usepackage{enumerate}
\usepackage{booktabs}
\usepackage[colorlinks=true,linkcolor=blue,citecolor=blue,urlcolor=blue]{hyperref}

\newtheorem{theorem}{Theorem}[section]
\newtheorem{lemma}[theorem]{Lemma}
\newtheorem{proposition}[theorem]{Proposition}
\newtheorem{corollary}[theorem]{Corollary}
\newtheorem{mainthm}{Theorem}

\theoremstyle{definition}
\newtheorem{definition}[theorem]{Definition}

\theoremstyle{remark}
\newtheorem{remark}[theorem]{Remark}

\numberwithin{equation}{section}

\newcommand{\FF}{\mathbb{F}}
\newcommand{\ZZ}{\mathbb{Z}}
\newcommand{\GL}{\mathrm{GL}}
\newcommand{\GammaL}{\Gamma\mathrm{L}}
\newcommand{\Sym}{\mathrm{Sym}}
\newcommand{\Stab}{\mathrm{Stab}}
\newcommand{\supp}{\mathrm{supp}}

\newcommand{\C}{\mathbf{C}}
\newcommand{\N}{\mathbf{N}}
\newcommand{\Z}{\mathbf{Z}}
\newcommand{\Pow}{\mathcal{P}}
\newcommand{\one}{\mathbf{1}}

\begin{document}

\title{Linear group actions with exactly one regular orbit}

\author[T. M. Keller]{Thomas Michael Keller}
\address{Department of Mathematics, Texas State University, 601 University Drive, San Marcos, TX 78666, USA}
\email{keller@txstate.edu}

\author[J. M. Riedl]{Jeffrey M. Riedl}
\address{Department of Mathematics, University of Akron, Akron, OH, USA}
\email{riedl@uakron.edu}



\begin{abstract}
Let $p$ be an odd prime, $G$ a finite $p$-group and $V$ a finite faithful irreducible $G$-module. We show that if $G$ has exactly one regular orbit on $V$, then $p=2^n-1$ is a Mersenne prime and $(G,V)$ embeds into $(P_m,V_m)$ for some $m\ge 1$, where $P_m$ is a Sylow $p$-subgroup of the monomial group $\GL_1(2^n)\wr S_{p^{m-1}}$ and $V_m$ is its natural module, regarded over $\FF_2$. We determine the largest orbit size $p^{x_m}$ of $P_m$ on $V_m$ and the number of orbits of this size, and we construct subgroups of $P_m$ of order $p^{x_m}$ whose unique regular orbit is a largest orbit of $P_m$; for $p=3$ such subgroups exist for every $m$, whereas for Mersenne primes $p\ge 7$ we prove that they exist only for $m\le 3$. Finally, we classify the subgroups of $P_m$ of order $p^{x_m}$ with exactly one regular orbit for $m=2$ and, with the help of a computer, for $p=3$ and $m=3$, where there are $13$ such pairs up to isomorphism; for $p=3$ and $m=4$ we show that there are exactly $4914$ conjugacy classes of such subgroups in $P_4$. In particular, such a subgroup is in general not determined by its regular orbit.
\end{abstract}

\maketitle

\section{Introduction}\label{sec:intro}

Orbit structures of finite linear groups, and in particular the existence of regular orbits, play an important role in the representation theory of finite groups; see, for instance, the books \cite{MW} and \cite{CY}. Usually one asks whether a group has a regular orbit on a given module. Many inductive regular-orbit arguments for imprimitive modules need the stabilizer of a block to have at least two regular orbits on its component which is why a certain focus has been to find more than one regular orbit. \\
In this paper, we initiate the study of the extreme situation in which there is exactly one regular orbit. More precisely, we are interested in pairs $(G,V)$, where $G$ is a finite $p$-group and $V$ is a finite faithful irreducible $G$-module, such that $G$ has exactly one regular orbit on $V$. The cases $p=2$ and $p$ odd are very different in nature, and in this paper we only treat odd primes~$p$.

Our first result shows that such pairs are rare and can be described in terms of iterated wreath products. Let $p=2^n-1$ be a Mersenne prime and $m\ge 1$, let $P_m\cong C_p\wr\cdots\wr C_p$ ($m$ factors) be a Sylow $p$-subgroup of the monomial group $\GL_1(2^n)\wr S_{p^{m-1}}\le \GL_{p^{m-1}}(2^n)$, and let $V_m$ be its natural module, regarded as a module over $\FF_2$ (see Section~\ref{sec:orbits} for details).

\begin{mainthm}\label{thmA}
Let $p$ be an odd prime, let $G\ne 1$ be a $p$-group, and let $V$ be a finite faithful irreducible $G$-module. If $G$ has exactly one regular orbit on $V$, then $p=2^n-1$ is a Mersenne prime, $V$ is a faithful irreducible $\FF_2G$-module, and $(G,V)$ is isomorphic to $(T,V_m)$ for some $m\ge 1$ and some subgroup $T\le P_m$.
\end{mainthm}

As a regular orbit of a subgroup $T\le P_m$ is contained in an orbit of $P_m$, in this situation $|T|$ is at most the largest orbit size $p^{x_m}$ of $P_m$ on $V_m$. Our second result determines this number, as well as the number of orbits of this size.

\begin{mainthm}\label{thmB}
Let $m\ge 2$. Then the largest orbit size of $P_m$ on $V_m$ is $p^{x_m}$, where
\[
x_m=p^{m-1}+\frac{p^{m-2}-1}{p-1}=\log_p|P_m|-p^{m-2}.
\]
The number $r_m$ of orbits of this size equals the number of regular orbits of $P_{m-2}$ on the power set of a set with $p^{m-2}$ elements; it satisfies $r_2=2$ and $r_{m+1}=(r_m^p-r_m)/p$. In particular, if $p=3$, then $P_m$ has exactly two orbits of size $3^{x_m}$ on $V_m$.
\end{mainthm}

We then concentrate on the extreme case of subgroups $T\le P_m$ of order $p^{x_m}$, so that the regular orbits of $T$ are largest orbits of $P_m$.

\begin{mainthm}\label{thmC}
Let $m\ge 2$. There are subgroups $T_U$ and $T_D$ of $P_m$ of order $p^{x_m}$, constructed explicitly in Section~\ref{sec:unique}, each of which has a largest orbit of $P_m$ on $V_m$ as a regular orbit; if $m\ge 3$, then $V_m$ is an irreducible module for both of them. If $p=3$, then each of $T_U$ and $T_D$ has exactly one regular orbit, and these are the two largest orbits of $P_m$ on $V_m$. Moreover, for every Mersenne prime $p$ and every $m\le 3$ there is a subgroup of $P_m$ of order $p^{x_m}$ that acts irreducibly on $V_m$ and has exactly one regular orbit.
\end{mainthm}

The main tool is a criterion (Theorem~\ref{thm:TDelta}) deciding when certain explicit subgroups $T_\Delta$ of $P_m$ have a given largest orbit of $P_m$ as a regular orbit. It is based on the uniserial structure of the permutation module of $P_k$ over $\FF_p$ and is closely related to the MDS property of certain repeated-root cyclic codes.

It is natural to ask whether a subgroup of $P_m$ of order $p^{x_m}$ with exactly one regular orbit is determined, up to conjugacy, by its regular orbit. This is true for $m\le 2$, but not in general.

\begin{mainthm}\label{thmD}
\begin{enumerate}[\rm(a)]
\item The subgroups of $P_2$ of order $p^{x_2}=p^p$ having exactly one regular orbit on $V_2$ are the base group of $P_2\cong C_p\wr C_p$ and the unique nonabelian maximal subgroup of $P_2$ of exponent $p$. Only the latter acts irreducibly on $V_2$.
\item Let $p=3$. The subgroups of $P_3$ of order $3^{x_3}=3^{10}$ having exactly one regular orbit on $V_3$ form $30$ conjugacy classes of subgroups of $P_3$. All of them act irreducibly on $V_3$, and they give rise to exactly $13$ isomorphism types of pairs $(T,V_3)$.
\item Let $p=3$. The subgroups of $P_4$ of order $3^{x_4}=3^{31}$ having exactly one regular orbit on $V_4$ form $4914$ conjugacy classes of subgroups of $P_4$, half of them with each of the two largest orbits of $P_4$ as regular orbit. All of them act irreducibly on $V_4$.
\item Let $p=3$ and $m\ge 3$. Then there are subgroups of $P_m$ of order $3^{x_m}$ that have the same unique regular orbit as $T_U$, act irreducibly on $V_m$, and are not conjugate to $T_U$ in $P_m$.
\end{enumerate}
\end{mainthm}

Parts~(b) and~(c) were obtained with the help of the computer algebra system GAP \cite{GAP}; see Sections~\ref{sec:m3} and~\ref{sec:m4}. For $m=4$ we have not determined the isomorphism types of the pairs $(T,V_4)$, and for $p=3$ the classification of all subgroups of $P_m$ of order $3^{x_m}$ with exactly one regular orbit remains open for $m\ge 5$.

Finally, we show that if some subgroup of $P_{m+1}$ has a largest orbit of $P_{m+1}$ as its unique regular orbit, then so does some subgroup of $P_m$ (Proposition~\ref{prop:descent}). Hence, for each Mersenne prime $p$, the values of $m$ for which such subgroups exist form an initial segment of the positive integers, which contains $\{1,2,3\}$ and, for $p=3$, all positive integers (Corollary~\ref{cor:descent}). Our last main result shows that for $p\ge 7$ this initial segment is exactly $\{1,2,3\}$.

\begin{mainthm}\label{thmE}
Let $p\ge 7$ be a Mersenne prime and $m\ge 4$. Then no subgroup of $P_m$ has a largest orbit of $P_m$ on $V_m$ as its unique regular orbit. Equivalently, every subgroup of $P_m$ of order $p^{x_m}$ has either no regular orbit or at least two regular orbits on $V_m$.
\end{mainthm}

Together with Theorem~\ref{thmC}, this shows that for a Mersenne prime $p\ge 7$ there is a subgroup of $P_m$ of order $p^{x_m}$ with exactly one regular orbit on $V_m$ if and only if $m\le 3$, and that in this case there is such a subgroup acting irreducibly on $V_m$. In the setting of Theorem~\ref{thmA}, if $p\ge 7$ and $|V|=2^{np^{m-1}}$ with $m\ge 4$, then $|G|<p^{x_m}$.

The proof of Theorem~\ref{thmE} in Section~\ref{sec:pge7} only uses linear algebra over $\FF_p$. We show that the regular orbits of a subgroup of $P_4$ of order $p^{x_4}$ are detected by the invertibility of certain matrices, which are circulant inside each block of $p$ points and constant between different blocks, and we then use a computation with small circulant matrices. The difference between $p=3$ and $p\ge 7$ comes from the number of orbits of $C_p$ on the proper nonempty subsets of a block of $p$ points: for $p=3$ there are only two such orbits, which is too few for this argument (Remark~\ref{rem:pge5}).

The paper is organized as follows. In Section~\ref{sec:orbits} we describe the orbits of $P_m$ on $V_m$ and prove Theorem~\ref{thmB}. In Section~\ref{sec:m2} we determine the regular orbits of all maximal subgroups of $P_2\cong C_p\wr C_p$, which settles the case $m=2$. In Section~\ref{sec:unique} we prove Theorem~\ref{thmA}, construct the groups $T_U$ and $T_D$, and prove Theorems~\ref{thmC} and~\ref{thmD}. Section~\ref{sec:pge7} contains the proof of Theorem~\ref{thmE}, and Section~\ref{sec:concl} collects some open problems.

\section{The largest orbits of \texorpdfstring{$P_m$}{Pm} on \texorpdfstring{$V_m$}{Vm}}\label{sec:orbits}

\subsection{Notation}
Throughout the paper, $p=2^n-1$ is a Mersenne prime (except in Lemma~\ref{lem:primefield} and at the beginning of Section~\ref{sec:unique}), $F=\FF_{2^n}$, and $\zeta$ is a generator of the multiplicative group $F^\times$, which has order~$p$. All groups act from the right. For a group $G$ acting on a set $\Omega$ we write $\Pow(\Omega)$ for the power set of $\Omega$, on which $G$ acts naturally, and $\Stab_G(\Gamma)$ for the set stabilizer of $\Gamma\subseteq\Omega$; a subset $\Gamma$ lies in a regular orbit of $G$ on $\Pow(\Omega)$ if and only if $\Stab_G(\Gamma)=1$. For a set $\Omega$ we denote by $\FF_p\Omega$ the $\FF_p$-vector space with basis $(e_\omega)_{\omega\in\Omega}$ and, for $\Gamma\subseteq\Omega$, by $\FF_p\Gamma$ the subspace spanned by the $e_\gamma$ with $\gamma\in\Gamma$; if a group acts on $\Omega$, then $\FF_p\Omega$ is the corresponding permutation module.

For $k\ge 0$ let $\Omega_k$ be a set with $p^k$ elements and let $P_k$ be a Sylow $p$-subgroup of $\Sym(\Omega_k)$. Thus $P_0=1$ and $P_k\cong C_p\wr\cdots\wr C_p$ ($k$ factors). For $m\ge 1$ let
\[
V_m=F^{\Omega_{m-1}}
\]
be the $F$-vector space with basis $(e_\omega)_{\omega\in\Omega_{m-1}}$. Let $H_{m-1}\le\GL(V_m)$ be the group of diagonal matrices (with respect to this basis) with entries in $F^\times$, and let $P_{m-1}$ act on $V_m$ by permuting the basis vectors. Then
\[
P_m=H_{m-1}\rtimes P_{m-1}
\]
is the group of monomial matrices whose nonzero entries lie in $F^\times$ and whose underlying permutation lies in $P_{m-1}$. Since $P_m\cong C_p\wr P_{m-1}$, this is consistent with the notation above. We call $V_m$ the \emph{natural $P_m$-module}; it has $F$-dimension $p^{m-1}$, and we always regard it as an $\FF_2P_m$-module. By Lemma~\ref{lem:irred} below, $V_m$ is faithful and irreducible.

We use the following well-known facts on the structure of $P_k$ and $P_m$ (see, for example, \cite{Kal} and \cite{GHH}).

\begin{enumerate}[(1)]
\item Iterating the decomposition above, we obtain $P_m=H_0H_1\cdots H_{m-1}$, where $H_i$ is an elementary abelian group of order $p^{p^i}$, the \emph{$i$-th layer} of $P_m$. For $0\le j<m$, the product $H_0\cdots H_j$ is a subgroup of $P_m$ isomorphic to $P_{j+1}$, and $H_{j}\cdots H_{m-1}$ is a normal subgroup of $P_m$. In particular $P_{m-1}=H_0\cdots H_{m-2}$.
\item Let $m\ge 2$. The orbits of $H_{m-2}$ on $\Omega_{m-1}$ form a system of $p^{m-2}$ blocks of size $p$ for $P_{m-1}$, the \emph{minimal blocks}, and $H_{m-2}$ is the kernel of the action of $P_{m-1}$ on the set of minimal blocks. We identify the set of minimal blocks with $\Omega_{m-2}$, so that $P_{m-1}$ acts on $\Omega_{m-2}$ as $P_{m-1}/H_{m-2}\cong P_{m-2}$, and we put
\[
Q=H_0\cdots H_{m-3}
\]
(so $Q=1$ if $m=2$); then $Q$ maps isomorphically onto $P_{m-1}/H_{m-2}$. For $i\in\Omega_{m-2}$ let $\beta_i\subseteq\Omega_{m-1}$ be the corresponding minimal block. Since $P_{m-1}\cong C_p\wr P_{m-2}$ with base group $H_{m-2}$, there are elements $y_i\in H_{m-2}$ ($i\in\Omega_{m-2}$), where $y_i$ acts as a $p$-cycle on $\beta_i$ and trivially on $\Omega_{m-1}\setminus\beta_i$, such that $H_{m-2}=\prod_i\langle y_i\rangle$ and $y_i^g=y_{ig}$ for all $g\in Q$.
\item Let $k\ge 1$. Then $\Omega_k$ is the disjoint union of $p$ blocks $\Omega_k^{(1)},\dots,\Omega_k^{(p)}$ of size $p^{k-1}$, the \emph{top blocks}, which are permuted cyclically by the top layer of $P_k$; the kernel $E$ of the action of $P_k$ on the top blocks is the direct product $E=E_1\times\dots\times E_p$, where $E_j$ acts on $\Omega_k^{(j)}$ as a Sylow $p$-subgroup of $\Sym(\Omega_k^{(j)})$ and trivially on the other top blocks. Thus $P_k\cong P_{k-1}\wr C_p$. We identify each pair $(\Omega_k^{(j)},E_j)$ with $(\Omega_{k-1},P_{k-1})$ by means of the top layer. For $k\ge 2$ every minimal block of $\Omega_k$ is contained in a top block, and this is compatible with the identifications above.
\item $P_k$ contains a $p^k$-cycle (indeed, any $p^k$-cycle generates a $p$-subgroup of $\Sym(\Omega_k)$, and all Sylow $p$-subgroups of $\Sym(\Omega_k)$ are conjugate).
\end{enumerate}

The map $\sum_\omega a_\omega e_\omega\mapsto \mathrm{diag}(\zeta^{a_\omega})_\omega$ is an isomorphism of $\FF_pP_{m-1}$-modules $\FF_p\Omega_{m-1}\to H_{m-1}$, and for $m\ge 2$ the map $\sum_i a_ie_i\mapsto\prod_i y_i^{a_i}$ is an isomorphism of $\FF_pP_{m-1}$-modules $\FF_p\Omega_{m-2}\to H_{m-2}$ (where $P_{m-1}$ acts on $H_{m-2}$ by conjugation and on $\FF_p\Omega_{m-2}$ through $P_{m-2}$). We use these identifications freely; in particular, we sometimes write $H_{m-1}$ and $H_{m-2}$ additively, and for $\Gamma\subseteq\Omega_{m-1}$ (resp.\ $\Gamma\subseteq\Omega_{m-2}$) we regard $\FF_p\Gamma$ as a subgroup of $H_{m-1}$ (resp.\ $H_{m-2}$).

For $v=\sum_\omega v_\omega e_\omega\in V_m$ the \emph{support} of $v$ is $\supp(v)=\{\omega\in\Omega_{m-1}: v_\omega\ne 0\}$. If $m\ge 2$ we also put
\[
\Delta(v)=\{\,i\in\Omega_{m-2}:\ \beta_i\subseteq\supp(v)\,\},
\]
the set of minimal blocks contained in the support of $v$.

\subsection{Orbits and supports}

\begin{lemma}\label{lem:supports}
Let $m\ge 1$, let $v\in V_m$ and put $S=\supp(v)$.
\begin{enumerate}[\rm(a)]
\item $H_{m-1}$ acts transitively on the set of vectors with support $S$, and $\C_{H_{m-1}}(v)=\FF_p(\Omega_{m-1}\setminus S)$.
\item Two vectors in $V_m$ lie in the same $P_m$-orbit if and only if their supports lie in the same $P_{m-1}$-orbit.
\item The image of $\C_{P_m}(v)$ in $P_m/H_{m-1}=P_{m-1}$ is $\Stab_{P_{m-1}}(S)$. In particular,
\[
|\C_{P_m}(v)|=p^{p^{m-1}-|S|}\,|\Stab_{P_{m-1}}(S)|\qquad\text{and}\qquad |v^{P_m}|=p^{|S|}\,|S^{P_{m-1}}|.
\]
\end{enumerate}
\end{lemma}

\begin{proof}
Part (a) holds because $F^\times$ acts regularly on $F\setminus\{0\}$. The elements of $P_m$ map supports to supports via their permutation parts. Hence, if $g\in P_{m-1}$ and $\supp(v')=Sg$, then $vg$ and $v'$ have the same support, so $vgh=v'$ for some $h\in H_{m-1}$ by (a); this proves (b), and (c) follows in the same way, applied to $v'=v$.
\end{proof}

For $k\ge 0$ let $s_k$ denote the number of regular orbits of $P_k$ on $\Pow(\Omega_k)$.

\begin{lemma}\label{lem:powerset}
\begin{enumerate}[\rm(a)]
\item Let $k\ge 0$ and $\Gamma\subseteq\Omega_{k+1}$, and put $\Gamma_j=\Gamma\cap\Omega_{k+1}^{(j)}$ for $1\le j\le p$. Then $\Gamma$ lies in a regular orbit of $P_{k+1}$ on $\Pow(\Omega_{k+1})$ if and only if each $\Gamma_j$, regarded as a subset of $\Omega_k$, lies in a regular orbit of $P_k$ on $\Pow(\Omega_k)$, and $\Gamma_1,\dots,\Gamma_p$ do not all lie in the same $P_k$-orbit.
\item $s_0=2$ and $s_{k+1}=(s_k^p-s_k)/p$ for $k\ge 0$. In particular $s_k\ge 2$ for all $k$.
\item If $\Gamma$ lies in a regular orbit of $P_k$ on $\Pow(\Omega_k)$, then so does $\Omega_k\setminus\Gamma$, and $\Gamma$ and $\Omega_k\setminus\Gamma$ lie in different $P_k$-orbits.
\item Every $p$-subgroup of $\Sym(\Omega)$, where $\Omega$ is a finite set, has a regular orbit on $\Pow(\Omega)$.
\end{enumerate}
\end{lemma}

\begin{proof}
(a) Let $E=E_1\times\dots\times E_p$ be the kernel of the action of $P_{k+1}$ on its top blocks. Then $\Stab_E(\Gamma)=\prod_j\Stab_{E_j}(\Gamma_j)$, so $\Stab_E(\Gamma)=1$ if and only if each $\Gamma_j$ lies in a regular orbit of $E_j\cong P_k$. Assume this is the case. Then the $E$-orbit of $\Gamma$ consists of all sets $\Gamma'$ such that $\Gamma'\cap\Omega_{k+1}^{(j)}$ lies in the $E_j$-orbit of $\Gamma_j$ for each $j$. If $g\in\Stab_{P_{k+1}}(\Gamma)\setminus E$, then $g$ permutes the top blocks as a $p$-cycle, and it follows that $\Gamma_1,\dots,\Gamma_p$ lie in the same $P_k$-orbit. Conversely, if they do, then a generator $c$ of the top layer maps $\Gamma$ to a set in the $E$-orbit of $\Gamma$, so that $\Gamma ce=\Gamma$ for some $e\in E$, and $ce\ne 1$.

(b) Clearly $s_0=2$. By (a), the sets lying in regular $P_{k+1}$-orbits form the union of the $E$-orbits $\mathcal{O}_1\times\dots\times\mathcal{O}_p$ (in the obvious sense), where the $\mathcal{O}_j$ are regular $P_k$-orbits, not all equal. There are $s_k^p-s_k$ such $E$-orbits, and the top layer permutes them in orbits of length~$p$. Thus $s_{k+1}=(s_k^p-s_k)/p$. As $(2^p-2)/p\ge 2$ and $(s^p-s)/p$ is increasing in $s\ge 1$, we obtain $s_k\ge 2$ by induction.

(c) We have $\Stab_{P_k}(\Gamma)=\Stab_{P_k}(\Omega_k\setminus\Gamma)$, and $|\Gamma|\ne|\Omega_k\setminus\Gamma|$ as $|\Omega_k|=p^k$ is odd.

(d) Let $G\le\Sym(\Omega)$ be a $p$-group and let $\Lambda_1,\dots,\Lambda_t$ be its orbits on $\Omega$. Each $\Lambda_i$ has size $p^{k_i}$ for some $k_i\ge 0$, and the image of $G$ in $\Sym(\Lambda_i)$ is contained in a Sylow $p$-subgroup of $\Sym(\Lambda_i)$, which is isomorphic (as a permutation group) to $P_{k_i}$. By (b), there is $\Gamma_i\subseteq\Lambda_i$ whose stabilizer in this Sylow subgroup is trivial. Then $\Gamma=\Gamma_1\cup\dots\cup\Gamma_t$ satisfies $\Stab_G(\Gamma)=1$.
\end{proof}

\begin{remark}
Part (d) of Lemma~\ref{lem:powerset} is a special case of a theorem of Gluck \cite[Corollary~1]{Gluck}, which states that every permutation group of odd order has a regular orbit on the power set.
\end{remark}

\begin{theorem}\label{thm:largest}
Let $m\ge 2$ and put
\[
x_m=\log_p|P_m|-p^{m-2}=p^{m-1}+\frac{p^{m-2}-1}{p-1}.
\]
\begin{enumerate}[\rm(a)]
\item Every orbit of $P_m$ on $V_m$ has size at most $p^{x_m}$.
\item Let $v\in V_m$. Then $|v^{P_m}|=p^{x_m}$ if and only if $|\supp(v)\cap\beta_i|\ge p-1$ for all $i\in\Omega_{m-2}$ and $\Delta(v)$ lies in a regular orbit of $P_{m-2}$ on $\Pow(\Omega_{m-2})$.
\item Two vectors $v,v'$ as in \textup{(b)} lie in the same $P_m$-orbit if and only if $\Delta(v)$ and $\Delta(v')$ lie in the same $P_{m-2}$-orbit. Hence the number of orbits of $P_m$ on $V_m$ of size $p^{x_m}$ is $r_m:=s_{m-2}$.
\item Let $\Delta\subseteq\Omega_{m-2}$ lie in a regular orbit of $P_{m-2}$. For each $i\in\Omega_{m-2}\setminus\Delta$ choose $\sigma_i\in\beta_i$, put
\[
S_\Delta=\bigcup_{i\in\Delta}\beta_i\ \cup\bigcup_{i\notin\Delta}(\beta_i\setminus\{\sigma_i\}),\qquad v_\Delta=\sum_{\omega\in S_\Delta}e_\omega,
\]
and let $z_i\in H_{m-1}$ be the diagonal matrix with entry $\zeta$ at $\sigma_i$ and entries $1$ elsewhere. Then
\[
\C_{P_m}(v_\Delta)=Y_\Delta\times Z_\Delta,\qquad\text{where } Y_\Delta=\prod_{i\in\Delta}\langle y_i\rangle\le H_{m-2}\ \text{ and }\ Z_\Delta=\prod_{i\notin\Delta}\langle z_i\rangle\le H_{m-1}.
\]
Every vector $v$ with $|v^{P_m}|=p^{x_m}$ and $\Delta(v)=\Delta$ is $P_m$-conjugate to $v_\Delta$; in particular $\C_{P_m}(v)$ is an elementary abelian subgroup of $H_{m-2}H_{m-1}$ of order $p^{p^{m-2}}$.
\end{enumerate}
\end{theorem}

\begin{proof}
Let $v\in V_m$ and $S=\supp(v)$. By Lemma~\ref{lem:supports}(c) we have $|v^{P_m}|=p^{x}$, where
\[
\log_p|P_m|-x=(p^{m-1}-|S|)+\log_p|\Stab_{P_{m-1}}(S)|.
\]
A nontrivial power of a $p$-cycle on $\beta_i$ stabilizes no proper nonempty subset of $\beta_i$. Hence $\Stab_{H_{m-2}}(S)=\prod_{i\in I_0}\langle y_i\rangle$, where $I_0$ is the set of all $i$ with $S\cap\beta_i\in\{\emptyset,\beta_i\}$. Therefore
\[
\log_p|P_m|-x\ \ge\ \sum_{i\in\Omega_{m-2}}\bigl(|\beta_i\setminus S|+[\,i\in I_0\,]\bigr)\ \ge\ p^{m-2},
\]
where $[\,i\in I_0\,]$ is $1$ if $i\in I_0$ and $0$ otherwise, because every summand is at least~$1$. Equality holds throughout if and only if $|S\cap\beta_i|\ge p-1$ for all $i$ and $\Stab_{P_{m-1}}(S)\le H_{m-2}$. Assume that $|S\cap\beta_i|\ge p-1$ for all $i$. Then $\Delta(v)$ is the set of minimal blocks contained in $S$, so the image of $\Stab_{P_{m-1}}(S)$ in $P_{m-1}/H_{m-2}\cong P_{m-2}$ stabilizes $\Delta(v)$. Conversely, let $g\in P_{m-1}$ be such that its image stabilizes $\Delta(v)$. Then $Sg$ contains the same minimal blocks as $S$ and meets each of the other minimal blocks in $p-1$ points. As $H_{m-2}$ acts transitively on the sets with this property, we have $Sgh=S$ for some $h\in H_{m-2}$. Hence the image of $\Stab_{P_{m-1}}(S)$ in $P_{m-2}$ is exactly $\Stab_{P_{m-2}}(\Delta(v))$, which proves (a) and (b); note that equality in (a) is attained because $s_{m-2}\ge 1$ by Lemma~\ref{lem:powerset}(b).

(c) By Lemma~\ref{lem:supports}(b), $v$ and $v'$ lie in the same $P_m$-orbit if and only if $\supp(v')=\supp(v)g$ for some $g\in P_{m-1}$. If this holds, then $\Delta(v')$ is the image of $\Delta(v)$ under $g$. Conversely, if $\Delta(v')$ is the image of $\Delta(v)$ under some $g\in P_{m-1}$, then $\supp(v)g$ and $\supp(v')$ contain the same minimal blocks and meet each of the other minimal blocks in $p-1$ points, so that $\supp(v')=\supp(v)gh$ for some $h\in H_{m-2}$.

(d) The elements $y_i$ ($i\in\Delta$) and $z_i$ ($i\notin\Delta$) fix $v_\Delta$, and they generate an elementary abelian group of order $p^{p^{m-2}}$, which by (b) and Lemma~\ref{lem:supports}(c) is the order of $\C_{P_m}(v_\Delta)$. The remaining assertions follow from (c) and Lemma~\ref{lem:supports}(b), and from the fact that $\C_{P_m}(v)$ is conjugate to $\C_{P_m}(v_\Delta)\le H_{m-2}H_{m-1}$, which is a normal subgroup of $P_m$.
\end{proof}

For $\Delta$ as in Theorem~\ref{thm:largest}(d) we denote by $O_\Delta=v_\Delta^{P_m}$ the corresponding largest orbit; it depends only on the $P_{m-2}$-orbit of $\Delta$, and $O_\Delta=\{v\in V_m:|v^{P_m}|=p^{x_m},\ \Delta(v)\in\Delta^{P_{m-2}}\}$.

\begin{corollary}\label{cor:recursion}
Put $x_1=1$ and $r_1=1$. Then $p^{x_m}$ is the largest orbit size of $P_m$ on $V_m$ and $r_m$ is the number of orbits of this size, for all $m\ge 1$. Moreover:
\begin{enumerate}[\rm(a)]
\item $x_2=p$ and $x_{m+1}=px_m+1$ for $m\ge 2$;
\item $r_2=2$ and $r_{m+1}=(r_m^p-r_m)/p$ for $m\ge 2$; in particular $r_m\ge 2$ for $m\ge 2$, and $r_m=2$ for all $m\ge 2$ if $p=3$;
\item for every $k\ge 0$, the group $P_k$ has exactly $r_{k+2}$ regular orbits on $\Pow(\Omega_k)$.
\end{enumerate}
\end{corollary}

\begin{proof}
For $m=1$ the group $P_1=F^\times$ acts regularly on $V_1\setminus\{0\}=F\setminus\{0\}$. For $m\ge 2$, everything follows from Theorem~\ref{thm:largest} and Lemma~\ref{lem:powerset}(b); for $p=3$ note that $(2^3-2)/3=2$.
\end{proof}

Note that the recursion $x_{m+1}=px_m+1$ fails for $m=1$, because $P_1$ has only one largest orbit on $V_1$ (that is, $r_1=1<2$).

\subsection{Vectors of types \texorpdfstring{$U_m$}{Um} and \texorpdfstring{$D_m$}{Dm}}

\begin{definition}\label{def:balanced}
A subset $X\subseteq\Omega_k$ is \emph{balanced} if either $k=0$, or $k\ge 1$ and the sets $X\cap\Omega_k^{(j)}$, $1\le j\le p$, have sizes differing by at most one and are balanced as subsets of $\Omega_k^{(j)}\cong\Omega_{k-1}$.
\end{definition}

Note that every subset of $\Omega_0$ and of $\Omega_1$ is balanced.

\begin{definition}\label{def:types}
Let $m\ge 2$. We define recursively when a vector $v\in V_m$ is of \emph{type $U_m$} or of \emph{type $D_m$}. If $m=2$, then $v$ is of type $U_2$ if $\supp(v)=\Omega_1$, and of type $D_2$ if $|\supp(v)|=p-1$. Now let $m\ge 3$. For $1\le j\le p$ let $X_j$ be the subspace of $V_m$ spanned by the $e_\omega$ with $\omega\in\Omega_{m-1}^{(j)}$. Then $V_m=X_1\oplus\dots\oplus X_p$, the top layer $H_0$ of $P_m$ permutes the $X_j$ cyclically, and the stabilizer $H_1\cdots H_{m-1}$ of all $X_j$ is a direct product of groups $E_1,\dots,E_p\cong P_{m-1}$ such that $X_j$ is the natural $E_j$-module. Write $v=v_1+\dots+v_p$ with $v_j\in X_j$. Then $v$ is of type $U_m$ (resp.\ $D_m$) if each $v_j$ is of type $U_{m-1}$ or $D_{m-1}$, exactly $(p+1)/2$ of them being of type $U_{m-1}$ (resp.\ $D_{m-1}$).
\end{definition}

The notation reminds us that for $m=2$ the centralizer in $P_2$ of a vector of type $U_2$ is not contained in the base group $H_1$ of $P_2$ (it is ``up''), whereas the centralizer of a vector of type $D_2$ is (it is ``down''); see Theorem~\ref{thm:largest}(d). Clearly the types are invariant under the action of $P_m$.

\begin{lemma}\label{lem:types}
Let $m\ge 2$ and let $v\in V_m$ be of type $U_m$ (resp.\ $D_m$). Then $v$ lies in a largest orbit of $P_m$ on $V_m$, the set $\Delta(v)$ is balanced, and $|\Delta(v)|=(p^{m-2}+1)/2$ (resp.\ $(p^{m-2}-1)/2$). If $p=3$, then the vectors of type $U_m$ form one of the two largest orbits of $P_m$ on $V_m$, denoted by $O_U$, and the vectors of type $D_m$ form the other one, denoted by $O_D$.
\end{lemma}

\begin{proof}
We argue by induction on $m$. For $m=2$ we have $\Omega_0=\{i\}$ and $\beta_i=\Omega_1$, so $\Delta(v)=\{i\}$ if $v$ is of type $U_2$ and $\Delta(v)=\emptyset$ if $v$ is of type $D_2$; by Theorem~\ref{thm:largest}(b), $v$ lies in a largest orbit. Let $m\ge 3$ and write $v=v_1+\dots+v_p$ as in Definition~\ref{def:types}. As the minimal blocks of $\Omega_{m-1}$ are contained in top blocks, $\Delta(v)$ is the disjoint union of the sets $\Delta(v_j)\subseteq\Omega_{m-2}^{(j)}$. By induction, each $\Delta(v_j)$ is balanced of size $(p^{m-3}\pm 1)/2$, every minimal block meets $\supp(v)$ in at least $p-1$ points, and each $\Delta(v_j)$ lies in a regular orbit of $P_{m-3}$. Hence $\Delta(v)$ is balanced, and its size is as claimed. Since the sets $\Delta(v_j)$ do not all have the same size, they do not all lie in the same $P_{m-3}$-orbit, so $\Delta(v)$ lies in a regular orbit of $P_{m-2}$ by Lemma~\ref{lem:powerset}(a). By Theorem~\ref{thm:largest}(b), $v$ lies in a largest orbit.

If $p=3$, then $P_m$ has exactly two largest orbits on $V_m$ by Corollary~\ref{cor:recursion}(b). Vectors of type $U_m$ and of type $D_m$ lie in largest orbits with $|\Delta(v)|=(3^{m-2}+1)/2$ and $|\Delta(v)|=(3^{m-2}-1)/2$, respectively, and hence in different ones. As the types are $P_m$-invariant, the last assertion follows.
\end{proof}

\subsection{An irreducibility criterion}

\begin{lemma}\label{lem:irred}
Let $m\ge 1$ and $T\le P_m$. Let $\bar T$ denote the image of $T$ in $P_m/H_{m-1}=P_{m-1}$ and put $M=T\cap H_{m-1}$. For $d\in H_{m-1}$ and $\omega\in\Omega_{m-1}$ let $d_\omega\in F^\times$ denote the diagonal entry of $d$ at $\omega$. Assume that
\begin{enumerate}[\rm(i)]
\item $\bar T$ acts transitively on $\Omega_{m-1}$, and
\item for every $\omega\in\Omega_{m-1}$ there is $d\in M$ with $d_\omega\ne 1$, and for any two distinct $\omega,\omega'\in\Omega_{m-1}$ there is $d\in M$ with $d_\omega\ne 1=d_{\omega'}$.
\end{enumerate}
Then $V_m$ is an irreducible $\FF_2T$-module. Condition \textup{(ii)} holds if $m\ge 2$ and $M$ contains the group of all $d\in H_{m-1}$ satisfying $\prod_{\omega\in\beta_i}d_\omega=1$ for all $i\in\Omega_{m-2}$. If \textup{(i)} does not hold, then $V_m$ is a reducible $\FF_2T$-module.
\end{lemma}

\begin{proof}
For $\omega\in\Omega_{m-1}$ the subspace $L_\omega=Fe_\omega$ is an $\FF_2M$-submodule of $V_m$, on which $d\in M$ acts as multiplication by $d_\omega\in F^\times$. By (ii), some element of $M$ acts on $L_\omega$ as multiplication by an element $\lambda\ne 1$ of $F^\times$, which is a primitive $p$-th root of unity. As the multiplicative order of $2$ modulo $p=2^n-1$ is $n$, we have $\FF_2(\lambda)=F$, and so $L_\omega$ is an irreducible $\FF_2M$-module. If $\omega\ne\omega'$, then by (ii) some $d\in M$ acts trivially on $L_{\omega'}$ but not on $L_\omega$, so $L_\omega\not\cong L_{\omega'}$. Hence $V_m=\bigoplus_\omega L_\omega$ is a direct sum of pairwise non-isomorphic irreducible $\FF_2M$-modules, and every $\FF_2M$-submodule of $V_m$ is the sum of some of the $L_\omega$. Since $T$ permutes the $L_\omega$ as $\bar T$ permutes $\Omega_{m-1}$, condition (i) implies that $V_m$ is irreducible for $T$.

If $M$ contains the group described in the lemma and $\omega\ne\omega'$, let $\omega''\in\beta_i\setminus\{\omega,\omega'\}$, where $\beta_i$ is the minimal block containing $\omega$ (this is possible as $p\ge 3$), and let $d$ have entry $\zeta$ at $\omega$, entry $\zeta^{-1}$ at $\omega''$ and entries $1$ elsewhere. Then $d\in M$ and $d_\omega\ne 1=d_{\omega'}$, so (ii) holds. Finally, if $\Lambda$ is an orbit of $\bar T$ on $\Omega_{m-1}$, then the $F$-span of the $e_\omega$ with $\omega\in\Lambda$ is $T$-invariant.
\end{proof}

In particular, $V_m$ is a faithful irreducible $\FF_2P_m$-module.

\section{Regular orbits of the maximal subgroups of \texorpdfstring{$C_p\wr C_p$}{Cp wr Cp}}\label{sec:m2}

In this section we treat the case $m=2$ in detail and determine the regular orbits of all maximal subgroups of $P_2\cong C_p\wr C_p$ on $V_2$. Throughout, we write $P=P_2$, $W=V_2$ and $B=H_1$. Thus $W=F^p$, and we write the elements of $W$ as $w=(v_1,\dots,v_p)$ with $v_j\in F$.

We regard $F$ as an $n$-dimensional $\FF_2$-space and let $x\in\GL(F)$ be the multiplication by $\zeta$. As $|F\setminus\{0\}|=p$, the group $\langle x\rangle$ acts regularly on $F^*=F\setminus\{0\}$. For $1\le j\le p$ let $x_j\in B$ act as $x$ on the $j$-th coordinate and trivially on the other coordinates. Then $B=\langle x_1,\dots,x_p\rangle$ is elementary abelian of order $p^p$. Let $y\in H_0$ be the permutation matrix mapping $(v_1,\dots,v_p)$ to $(v_p,v_1,\dots,v_{p-1})$. Then $y$ permutes $x_1,\dots,x_p$ cyclically by conjugation, $H_0=\langle y\rangle$, and $P=B\rtimes\langle y\rangle$ is the regular wreath product of order $p^{p+1}$.

Let $\ZZ_p=\{0,1,\dots,p-1\}$ be the group of integers modulo $p$, and let $\ZZ_p^p$ be the direct sum of $p$ copies of $\ZZ_p$. For $a=(a_1,\dots,a_p)\in\ZZ_p^p$ put
\[
\sigma(a)=a_1+\dots+a_p\in\ZZ_p\qquad\text{and}\qquad\eta(a)=x_1^{a_1}\cdots x_p^{a_p}\in B.
\]
Then $\sigma\colon\ZZ_p^p\to\ZZ_p$ is a homomorphism and $\eta\colon\ZZ_p^p\to B$ is an isomorphism. For $r\in\ZZ_p$ let $(\ZZ_p^p)_r=\sigma^{-1}(r)$; then $\ZZ_p^p$ is the disjoint union of the sets $(\ZZ_p^p)_r$, $r\in\ZZ_p$. Finally, let
\[
N=\eta\bigl((\ZZ_p^p)_0\bigr),\qquad z=x_1x_2\cdots x_p,\qquad\text{and}\qquad M_k=\langle N,\,yx_1^k\rangle\quad(k\in\ZZ_p).
\]

\begin{lemma}\label{lem:Pstructure}
\begin{enumerate}[\rm(a)]
\item $N$ is a normal subgroup of $P$ with $|B:N|=p$. Moreover, $N=[B,y]=P'=\Phi(P)$, and $\Z(P)=\langle z\rangle$.
\item Let $k\in\{1,\dots,p-1\}$ and $a\in\ZZ_p^p$. Then $(y^k\eta(a))^p=z^{\sigma(a)}$. Hence $y^k\eta(a)$ has order $p$ if $\sigma(a)=0$, and order $p^2$ otherwise.
\item The maximal subgroups of $P$ are $B$ and $M_0,M_1,\dots,M_{p-1}$. The group $B$ is elementary abelian, $M_0$ is nonabelian of exponent $p$, and $M_k$ has exponent $p^2$ for $k\ne 0$.
\end{enumerate}
\end{lemma}

\begin{proof}
(a) As $\sigma$ is invariant under permutations of the coordinates, $N$ is a $y$-invariant subgroup of index $p$ in $B$, hence normal in $P$. For $a\in\ZZ_p^p$ we have $[\eta(a),y]=\eta(a^y-a)$, where $a^y$ denotes the cyclic shift of $a$ corresponding to $y$. The kernel of $a\mapsto a^y-a$ consists of the constant vectors, so $[B,y]$ has order $p^{p-1}$; since $\sigma(a^y-a)=0$, we conclude that $[B,y]=N$. The images of $x_1$ and $y$ generate $P/N$ and commute, and both have order dividing $p$. Hence $P/N$ is elementary abelian of order $p^2$, which yields $\Phi(P)\le N=[B,y]\le P'\le\Phi(P)$. Finally, $\C_B(y)=\langle z\rangle$, and no element of $P\setminus B$ centralizes $B$, so $\Z(P)=\langle z\rangle$.

(b) For elements $h,b$ of any group and $l\ge 1$ we have
\[
(hb)^l=h^l\,(h^{-(l-1)}bh^{l-1})\cdots(h^{-2}bh^{2})(h^{-1}bh)\,b.
\]
We apply this with $h=y^k$, $b=\eta(a)$ and $l=p$. Then $h^p=1$, and since $y^k$ generates $\langle y\rangle$, the $p$ conjugates $h^{-i}bh^{i}$ are the elements $\eta(a')$, where $a'$ runs through the cyclic shifts of $a$. Their product is $\eta(\sigma(a),\dots,\sigma(a))=z^{\sigma(a)}$. As $z$ has order $p$, the assertion follows.

(c) The maximal subgroups of $P$ contain $\Phi(P)=N$ and correspond to the $p+1$ subgroups of order $p$ of $P/N\cong C_p\times C_p$, which are generated by the images of $x_1$ and of $yx_1^k$, $k\in\ZZ_p$. This yields $B,M_0,\dots,M_{p-1}$. By the identity in (b), $(yx_1^k)^j=y^j\eta(a)$ with $\sigma(a)=jk$; hence every element of $M_k\setminus N$ has the form $y^j\eta(a)$ with $j\ne 0$ and $\sigma(a)=jk$. Thus, by (b), $M_0$ has exponent $p$, and $M_k$ has exponent $p^2$ for $k\ne 0$. Finally $M_0$ is nonabelian, because $\C_N(y)=\langle z\rangle\ne N$.
\end{proof}

Let $W^*=W\setminus\{0\}$. For $1\le j\le p$ let $\mathcal{S}_j$ be the set of all $j$-element subsets of $\{1,\dots,p\}$, and put $A_j=\{w\in W^*:\ |\supp(w)|=j\}$ and, for $S\in\mathcal{S}_j$, $A_j(S)=\{w\in W^*:\ \supp(w)=S\}$. Then $W^*$ is the disjoint union of $A_1,\dots,A_p$, and $A_j$ is the disjoint union of the sets $A_j(S)$, $S\in\mathcal{S}_j$. We have $|A_j(S)|=p^j$ and $|A_j|=\binom pj p^j$. The sets $A_j$ are $P$-invariant, and the sets $A_j(S)$ are $B$-invariant, hence $N$-invariant.

\begin{lemma}\label{lem:Aj}
Let $1\le j\le p$ and $S\in\mathcal{S}_j$.
\begin{enumerate}[\rm(a)]
\item $B$ acts transitively on $A_j(S)$. In particular, $B$ acts regularly on $A_p$.
\item If $j<p$, then $A_j(S)$ is an $N$-orbit.
\item If $j<p$, $w\in A_j$ and $g\in P\setminus B$, then $\supp(w^g)\ne\supp(w)$.
\end{enumerate}
\end{lemma}

\begin{proof}
(a) The first assertion holds because $\langle x\rangle$ acts transitively on $F^*$. The second one follows, since $A_p=A_p(\{1,\dots,p\})$ and $|B|=p^p=|A_p|$.

(b) Let $w=(v_1,\dots,v_p)$ and $w'=(v_1',\dots,v_p')$ be elements of $A_j(S)$. By (a) there is $a=(a_1,\dots,a_p)\in\ZZ_p^p$ with $w^{\eta(a)}=w'$. Choose $i\notin S$, so that $v_i=0=v_i'$, and replace $a_i$ by the unique element of $\ZZ_p$ for which $\sigma(a)=0$. For the new $a$ we still have $w^{\eta(a)}=w'$, and now $\eta(a)\in N$.

(c) Write $g=\eta(a)y^k$ with $a\in\ZZ_p^p$ and $k\in\{1,\dots,p-1\}$. Since $w^{\eta(a)}=(v_1x^{a_1},\dots,v_px^{a_p})$ for $w=(v_1,\dots,v_p)$, we have $\supp(w^{\eta(a)})=\supp(w)$. Now $w^g$ arises from $w^{\eta(a)}$ by applying the $p$-cycle $y^k$ to the coordinates. As $w^{\eta(a)}$ has at least one zero and at least one nonzero coordinate, $\supp(w^g)\ne\supp(w^{\eta(a)})=\supp(w)$.
\end{proof}

\begin{lemma}\label{lem:Ajorbits}
Let $M\in\{M_0,\dots,M_{p-1}\}$ and $1\le j\le p-1$. Then every $M$-orbit on $A_j$ has size $p^{j+1}$. In particular, $M$ has no regular orbit on $A_j$ if $j\le p-2$, and $A_{p-1}$ is a regular orbit of $M$.
\end{lemma}

\begin{proof}
Let $\mathcal{O}$ be an $M$-orbit on $A_j$. As $N\le M$, the set $\mathcal{O}$ is a union of $N$-orbits, which by Lemma~\ref{lem:Aj}(b) are sets of the form $A_j(S)$, $S\in\mathcal{S}_j$. These are permuted transitively by $M/N$, which has order $p$. Hence $\mathcal{O}$ is the union of either one or $p$ such sets. Let $w\in\mathcal{O}$ and $g\in M\setminus B$. By Lemma~\ref{lem:Aj}(c), $w^g\in\mathcal{O}\setminus A_j(\supp(w))$. Therefore $|\mathcal{O}|=p\cdot p^j$. Since $|M|=p^p$, the orbit $\mathcal{O}$ is regular if and only if $j=p-1$, and since $|A_{p-1}|=p\cdot p^{p-1}$, the set $A_{p-1}$ is a single $M$-orbit.
\end{proof}

We now label the elements of $A_p$. Fix $u_0\in F^*$ and put $u_i=u_0x^i$ for $i\in\ZZ_p$. Then $u_0,\dots,u_{p-1}$ are the distinct elements of $F^*$, and $u_jx^i=u_{i+j}$ for $i,j\in\ZZ_p$. For $c=(c_1,\dots,c_p)\in\ZZ_p^p$ let
\[
w(c)=(u_{c_1},\dots,u_{c_p})\in A_p.
\]
The map $c\mapsto w(c)$ is a bijection $\ZZ_p^p\to A_p$. For $r\in\ZZ_p$ let
\[
E_r=\{\,w(c):\ c\in(\ZZ_p^p)_r\,\}.
\]
Then $|E_r|=p^{p-1}$, and $A_p$ is the disjoint union of $E_0,\dots,E_{p-1}$.

\begin{lemma}\label{lem:Er}
\begin{enumerate}[\rm(a)]
\item For $a,c\in\ZZ_p^p$ we have $w(c)^{\eta(a)}=w(a+c)$, and $w(c)^y=w(c^y)$, where $c^y$ is the corresponding cyclic shift of $c$. In particular, $E_r^{\eta(a)}=E_{r+\sigma(a)}$ and $E_r^y=E_r$ for all $r\in\ZZ_p$.
\item The orbits of $N$ on $A_p$ are $E_0,E_1,\dots,E_{p-1}$.
\end{enumerate}
\end{lemma}

\begin{proof}
(a) We have $w(c)^{\eta(a)}=(u_{c_1}x^{a_1},\dots,u_{c_p}x^{a_p})=(u_{a_1+c_1},\dots,u_{a_p+c_p})=w(a+c)$. The second formula is clear, and the remaining assertions follow as $\sigma$ is a homomorphism invariant under permutations of the coordinates.

(b) By (a), $N$ leaves each $E_r$ invariant. As $B$ acts regularly on $A_p$ by Lemma~\ref{lem:Aj}(a), $N$ acts semiregularly on $A_p$, so its orbits have size $p^{p-1}=|E_r|$.
\end{proof}

\begin{theorem}\label{thm:maxsubgroups}
\begin{enumerate}[\rm(a)]
\item $B$ has exactly one regular orbit on $W$, namely $A_p$.
\item $M_0$ has exactly one regular orbit on $W$, namely $A_{p-1}$. The orbits of $M_0$ on $A_p$ are $E_0,E_1,\dots,E_{p-1}$.
\item If $k\ne 0$, then $M_k$ has exactly two regular orbits on $W$, namely $A_{p-1}$ and $A_p$.
\item $W$ is an irreducible $M_k$-module for every $k\in\ZZ_p$, but it is a reducible $B$-module.
\end{enumerate}
\end{theorem}

\begin{proof}
(a) The $B$-orbits on $A_j$, $j<p$, are the sets $A_j(S)$, which have size $p^j<p^p=|B|$, and $B$ acts regularly on $A_p$ by Lemma~\ref{lem:Aj}(a).

(b) and (c) By Lemma~\ref{lem:Ajorbits}, the only regular orbit of $M_k$ contained in $A_1\cup\dots\cup A_{p-1}$ is $A_{p-1}$. As $M_0=\langle N,y\rangle$, Lemma~\ref{lem:Er} shows that the $M_0$-orbits on $A_p$ are the sets $E_r$, which have size $p^{p-1}$. Now let $k\ne 0$. The group $M_k/N$ of order $p$ permutes the $N$-orbits $E_0,\dots,E_{p-1}$, and by Lemma~\ref{lem:Er}(a) the element $yx_1^k=y\eta(k,0,\dots,0)$ maps $E_r$ to $E_{r+k}\ne E_r$. Hence $M_k$ acts transitively on $A_p$, and since $|M_k|=p^p=|A_p|$, it acts regularly.

(d) Let $k\in\ZZ_p$. As $|M_k\cap B|=p^{p-1}$ and $N\le M_k\cap B$, we have $M_k\cap B=N$, which is the group of all $d\in B$ with $\prod_jd_j=1$. Moreover, the image of $M_k$ in $P/B\cong\langle y\rangle$ is all of $P/B$, which acts transitively on the coordinates. Hence $W$ is irreducible for $M_k$ by Lemma~\ref{lem:irred}. Clearly $Fe_1$ is a $B$-invariant subspace of $W$.
\end{proof}

Since $x_2=p$ by Corollary~\ref{cor:recursion}, the subgroups of $P_2$ of order $p^{x_2}$ are precisely its maximal subgroups. Hence Theorem~\ref{thm:maxsubgroups} implies:

\begin{corollary}\label{cor:m2}
The subgroups of $P_2$ of order $p^{x_2}=p^p$ having exactly one regular orbit on $V_2$ are precisely $B=H_1$, with regular orbit $A_p$, and $M_0$, with regular orbit $A_{p-1}$. Only $M_0$ acts irreducibly on $V_2$.
\end{corollary}

Note that $A_p$ and $A_{p-1}$ are the two largest orbits of $P_2$ on $V_2$, consisting of the vectors of type $U_2$ and $D_2$, respectively.

\begin{remark}
The groups $M_1,\dots,M_{p-1}$ are conjugate in $\GL(W)$, namely by suitable permutation matrices normalizing $\langle y\rangle$, but they are not conjugate in $P$, as they are normal in $P$.
\end{remark}

\section{\texorpdfstring{$p$}{p}-groups with a unique regular orbit}\label{sec:unique}

\subsection{Reduction to subgroups of \texorpdfstring{$P_m$}{Pm}}

In this subsection, $p$ denotes an arbitrary odd prime. We first show that it suffices to consider modules over prime fields.

\begin{lemma}\label{lem:primefield}
Let $G\ne 1$ be a group of odd order, let $\FF$ be a finite field with prime field $\FF_q$, and let $V$ be a faithful irreducible $\FF G$-module. If $G$ has exactly one regular orbit on $V$, then $V$ is irreducible as an $\FF_qG$-module.
\end{lemma}

\begin{proof}
Let $W$ be an irreducible $\FF_qG$-submodule of $V$. For $\lambda\in\FF^\times$, the subspace $\lambda W$ is an $\FF_qG$-submodule isomorphic to $W$, and $\sum_\lambda\lambda W$ is a nonzero $\FF G$-submodule, hence equal to $V$. Thus $V\cong W^{e}$ as $\FF_qG$-modules for some $e\ge 1$, and $W$ is faithful. Suppose that $e\ge 2$, and write the elements of $V$ as $e$-tuples $(w_1,\dots,w_e)$ with $w_i\in W$, so that $G$ acts diagonally. If $(w_1,w_2,\dots,w_e)$ lies in a regular orbit, then so does $(w_2,w_1,w_3,\dots,w_e)$, as both vectors have the centralizer $\bigcap_i\C_G(w_i)$. If these vectors lie in the same orbit, then some $g\in G$ interchanges $w_1$ and $w_2$ and fixes $w_3,\dots,w_e$; then $g^2=1$, so $g=1$ as $|G|$ is odd, and $w_1=w_2$. Hence, if $G$ has exactly one regular orbit, then any vector $(w_1,\dots,w_e)$ in this orbit satisfies $w_1=w_2=\dots=w_e=:w$ (by the same argument for the other pairs of coordinates), and $\C_G(w)=1$. But then $(w,0,\dots,0)$ also lies in a regular orbit, and its coordinates are not all equal, since $w\ne 0$ as $G\ne 1$. This contradiction shows that $e=1$.
\end{proof}

\begin{theorem}\label{thm:reduction}
Let $p$ be an odd prime, let $G\ne 1$ be a $p$-group, and let $V$ be a finite faithful irreducible $G$-module such that $G$ has exactly one regular orbit on $V$. Then $p=2^n-1$ is a Mersenne prime, $V$ is a faithful irreducible $\FF_2G$-module, and there are $m\ge 1$ and an $\FF_2$-linear isomorphism $V\to V_m$ conjugating $G$ into $P_m$. In particular, $|G|\le p^{x_m}$.
\end{theorem}

\begin{proof}
By Lemma~\ref{lem:primefield} we may assume that $V$ is an irreducible $\FF_qG$-module for a prime $q$. As the only irreducible $\FF_pG$-module is the trivial one, $q\ne p$.

First assume that $V$ is a primitive $\FF_qG$-module. Then $V$ is quasi-primitive, and so every abelian normal subgroup of $G$ is cyclic (see the beginning of \S1 of \cite{MW}). As $p$ is odd, $G$ is cyclic by \cite[Corollary~1.3]{MW}. Let $|V|=q^d$. By \cite[Theorem~2.1]{MW}, applied with $A=G$, we may identify $V$ with $\FF_{q^d}$ in such a way that $G\le\Gamma_0(q^d)$, that is, $G$ acts by multiplication with elements of $\FF_{q^d}^\times$. In particular, $G$ acts semiregularly on $V\setminus\{0\}$, so the number of regular orbits of $G$ on $V$ is $(q^d-1)/|G|$, and therefore $q^d-1=|G|=p^a$ for some $a\ge 1$. Since $p$ is odd, \cite[Proposition~3.1]{MW} yields $q=2$, $a=1$, and that $p=2^d-1$ is a Mersenne prime. Thus $d=n$, $|G|=p$ and $(G,V)\cong(F^\times,F)=(P_1,V_1)$.

Now assume that $V$ is imprimitive, and choose a decomposition $V=V_1\oplus\dots\oplus V_k$, $k\ge 2$, into subspaces permuted transitively by $G$, such that $\dim V_1$ is minimal. Let $H$ be the stabilizer of $V_1$ in $G$. Then $V$ is induced from the $H$-module $V_1$ (see the beginning of \S1 of \cite{MW}), so that $V_1$ is a faithful irreducible module for $\bar H=H/\C_H(V_1)$, which is primitive by the minimality of $\dim V_1$, and $k=|G:H|$ is a power of $p$. Note that $\bar H\ne 1$, since otherwise $V_1$ would be the trivial module of dimension one and $V$ would be the permutation module of $G$ on $\{V_1,\dots,V_k\}$, which is reducible as $k\ge 2$. By the first part of the proof, $\bar H$ is cyclic and acts semiregularly on $V_1\setminus\{0\}$; in particular, $\bar H$ has at least one regular orbit on $V_1$.

We claim that $\bar H$ has exactly one regular orbit on $V_1$. Suppose that $u,u'\in V_1$ lie in different regular $\bar H$-orbits. Let $\Lambda=\{V_1,\dots,V_k\}$. By Lemma~\ref{lem:powerset}(d), the image of $G$ in $\Sym(\Lambda)$ has a regular orbit on $\Pow(\Lambda)$; let $\Gamma\subseteq\Lambda$ lie in such an orbit. Choose $g_1=1,g_2,\dots,g_k\in G$ with $V_i=V_1g_i$, and put
\[
v=\sum_{V_i\in\Gamma}ug_i+\sum_{V_i\notin\Gamma}u'g_i\qquad\text{and}\qquad v'=\sum_{V_i\in\Gamma}u'g_i+\sum_{V_i\notin\Gamma}ug_i.
\]
Let $g\in\C_G(v)$, and let $V_ig=V_j$. Then $h=g_igg_j^{-1}\in H$, and $g$ maps the component $ug_i$ (resp.\ $u'g_i$) of $v$ in $V_i$ to $(uh)g_j$ (resp.\ $(u'h)g_j$). As $u$ and $u'$ lie in different $H$-orbits, $g$ maps $\Gamma$ to $\Gamma$, so $g$ stabilizes every $V_i$. But then $g_igg_i^{-1}\in H$ fixes $u$ or $u'$ and therefore acts trivially on $V_1$, that is, $g$ acts trivially on $V_i$, for every $i$. Thus $g=1$, and $v$ lies in a regular orbit of $G$; similarly, so does $v'$. If $v'=vg$ for some $g\in G$, then the same argument shows that $g$ maps $\Gamma$ to $\Lambda\setminus\Gamma$, which is impossible as $|\Gamma|\ne k-|\Gamma|$ ($k$ is odd). Hence $G$ has at least two regular orbits on $V$, a contradiction.

By the first part of the proof, $q=2$, $p=2^n-1$ is a Mersenne prime, $|\bar H|=p$, and we may identify $V_1$ with $F$ in such a way that $\bar H=F^\times$ acts by multiplication. Let $S$ be the image of $G$ in $\Sym(\Lambda)\cong\Sym(k)$. By \cite[Lemma~2.8]{MW}, $G$ is isomorphic as a linear group to a subgroup of $\bar H\wr S$, acting on $V_1^k=F^k$; that is, $V$ can be identified with $F^k$ in such a way that $G$ consists of monomial matrices with nonzero entries in $F^\times$ and underlying permutations in $S$. As $S$ is a $p$-group, it is contained in a Sylow $p$-subgroup of $\Sym(k)$. Writing $k=p^{m-1}$ and identifying this Sylow subgroup with $P_{m-1}$, we obtain an $\FF_2$-linear isomorphism $V\to V_m$ conjugating $G$ into $P_m$. The final assertion holds because a regular orbit of a subgroup of $P_m$ is contained in an orbit of $P_m$.
\end{proof}

Theorem~\ref{thm:reduction} implies Theorem~\ref{thmA}. From now on, $p=2^n-1$ is again a Mersenne prime.

\begin{remark}\label{rem:sylow}
The group $P_m$ is a Sylow $p$-subgroup of $\GL(V_m)=\GL_{np^{m-1}}(2)$; compare \cite{Weir}. Indeed, $p$ divides $2^i-1$ if and only if $n$ divides $i$, and the $p$-part of $2^{nj}-1$ equals $p$ times the $p$-part of $j$, as $p$ is odd. Hence the $p$-part of $|\GL_{np^{m-1}}(2)|$ is $p^{p^{m-1}}$ times the $p$-part of $(p^{m-1})!$, which is $|P_m|$. Consequently, by Theorem~\ref{thm:reduction}, the isomorphism types of pairs $(G,V)$ as in Theorem~\ref{thmA} with $|V|=2^{np^{m-1}}$ correspond to the $\GL(V_m)$-conjugacy classes of subgroups of $P_m$ that act irreducibly on $V_m$ and have exactly one regular orbit on $V_m$. (Irreducibility is not automatic; see Corollary~\ref{cor:m2}.) Note that subgroups of $P_m$ may be conjugate in $\GL(V_m)$ without being conjugate in $P_m$.
\end{remark}

Let $m\ge 1$. By Theorem~\ref{thm:reduction}, a subgroup of $P_m$ with exactly one regular orbit has order at most $p^{x_m}$, and we now concentrate on the extreme case of subgroups of order $p^{x_m}$.

\begin{lemma}\label{lem:regular}
Let $T\le P_m$ with $|T|=p^{x_m}$. Then the regular orbits of $T$ on $V_m$ are exactly the largest orbits $O$ of $P_m$ on $V_m$ such that $T\cap\C_{P_m}(v)=1$ for some (equivalently, every) $v\in O$.
\end{lemma}

\begin{proof}
A regular orbit of $T$ has size $p^{x_m}$ and is contained in an orbit of $P_m$, which has size at most $p^{x_m}$; so it is a largest orbit of $P_m$. Conversely, let $O$ be a largest orbit of $P_m$ and $v\in O$. Then $v^T=O$ if and only if $|T:T\cap\C_{P_m}(v)|=p^{x_m}$, that is, if and only if $T\cap\C_{P_m}(v)=1$; and if $T$ acts transitively on $O$, then all point stabilizers in $T$ are conjugate.
\end{proof}

For $m\ge 2$ and a largest orbit $O$ of $P_m$ on $V_m$, let $\mathcal{S}(O)$ denote the set of supports of the vectors in $O$. By Lemma~\ref{lem:supports} and Theorem~\ref{thm:largest}, $\mathcal{S}(O_\Delta)$ is a single $P_{m-1}$-orbit of length $|P_{m-1}|/p^{|\Delta|}$.

The following lemma restricts the possible images in $P_{m-1}$ of subgroups of $P_m$ that act transitively on a largest orbit; we use the notation $J^c(\Omega_{m-2})$ of Lemma~\ref{lem:uniserial} below.

\begin{lemma}\label{lem:image}
Let $m\ge 2$, let $T\le P_m$ act transitively on a largest orbit $O_\Delta$ of $P_m$ on $V_m$, and let $\bar T$ be the image of $T$ in $P_m/H_{m-1}=P_{m-1}$. Then $\bar TH_{m-2}=P_{m-1}$, that is, $T$ maps onto $P_m/H_{m-2}H_{m-1}\cong P_{m-2}$, and $\bar T\cap H_{m-2}=J^c(\Omega_{m-2})$ for some $c\ge 0$. If moreover $|T|=p^{x_m}$, so that $O_\Delta$ is a regular orbit of $T$, then $c\le|\Delta|$.
\end{lemma}

\begin{proof}
The map $v\mapsto\Delta(v)$ maps $O_\Delta$ onto the $P_{m-2}$-orbit of $\Delta$ and is compatible with the natural homomorphism $P_m\to P_{m-2}$. As $T$ acts transitively on $O_\Delta$, its image in $P_{m-2}$ acts transitively on the orbit of $\Delta$, which is a regular $P_{m-2}$-orbit; hence this image is all of $P_{m-2}$, and $\bar TH_{m-2}=P_{m-1}$. As $H_{m-2}$ is abelian and normalized by $\bar T$, the group $\bar T\cap H_{m-2}$ is normal in $\bar TH_{m-2}=P_{m-1}$, so it equals $J^c(\Omega_{m-2})$ for some $c$ by Lemma~\ref{lem:uniserial}. Finally, $\bar T$ acts transitively on $\mathcal{S}(O_\Delta)$, so $|P_{m-2}|\,p^{p^{m-2}-c}=|\bar T|\ge|P_{m-1}|/p^{|\Delta|}=|P_{m-2}|\,p^{p^{m-2}-|\Delta|}$, that is, $c\le|\Delta|$.
\end{proof}

\subsection{Uniserial permutation modules and good sets}

\begin{lemma}\label{lem:uniserial}
Let $k\ge 0$. Then the $\FF_pP_k$-module $\FF_p\Omega_k$ is uniserial. More precisely, for $0\le d\le p^k$ it has a unique submodule $J^d=J^d(\Omega_k)$ of codimension $d$, we have $J^{d+1}=[J^d,P_k]$ for $d<p^k$, and $J^{p^k-1}$ is spanned by $\sum_{\omega\in\Omega_k}e_\omega$. If $g\in P_k$ is a $p^k$-cycle, then $J^d=\FF_p\Omega_k\cdot(g-1)^d$.
\end{lemma}

\begin{proof}
Let $g\in P_k$ be a $p^k$-cycle, and fix $\omega_0\in\Omega_k$. Then $\omega_0g^x\mapsto t^x$ ($0\le x<p^k$) defines an isomorphism of $\FF_p\langle g\rangle$-modules from $\FF_p\Omega_k$ onto $R=\FF_p[t]/(t^{p^k}-1)=\FF_p[t]/(t-1)^{p^k}$, on which $g$ acts as multiplication by $t$. The $\langle g\rangle$-submodules of $R$ are the ideals $(t-1)^dR$, $0\le d\le p^k$, which form a chain. Hence the $P_k$-submodules of $\FF_p\Omega_k$ form a chain as well. As $P_k$ is a $p$-group, $\FF_p\Omega_k$ has a composition series of $P_k$-submodules with one-dimensional factors, so there is a $P_k$-submodule of every codimension $d$, and it must coincide with $\FF_p\Omega_k\cdot(g-1)^d$. For $d<p^k$, the group $[J^d,P_k]$ is the smallest submodule $X$ of $J^d$ such that $P_k$ acts trivially on $J^d/X$; by uniseriality, this is the unique maximal submodule $J^{d+1}$ of $J^d$. Finally, the unique minimal submodule $J^{p^k-1}$ consists of the fixed points of $P_k$, which are spanned by $\sum_\omega e_\omega$ as $P_k$ is transitive.
\end{proof}

Lemma~\ref{lem:uniserial} is well known; it implies, for example, that for each divisor $p^i$ of $|H_{m-1}|$ there is exactly one subgroup of $H_{m-1}$ of order $p^i$ that is normal in $P_m$, namely $J^{p^{m-1}-i}(\Omega_{m-1})$.

\begin{lemma}\label{lem:MDS}
\begin{enumerate}[\rm(a)]
\item Let $f=\sum_{e\in E}c_et^e\in\FF_p[t]$, where $\emptyset\ne E\subseteq\{0,1,\dots,p-1\}$ and $c_e\ne 0$ for all $e\in E$. Then $(t-1)^{|E|}$ does not divide $f$.
\item Let $X\subseteq\Omega_1$. Then $\FF_p\Omega_1=J^{|X|}(\Omega_1)\oplus\FF_pX$.
\end{enumerate}
\end{lemma}

\begin{proof}
(a) We argue by induction on $|E|$. If $|E|=1$, then $f=c_et^e$ and $f(1)=c_e\ne 0$. Let $|E|>1$ and $e_0=\min E$. Then $f=t^{e_0}f_0$ with $f_0=\sum_{e\in E}c_et^{e-e_0}$, and $(t-1)^{|E|}$ divides $f$ if and only if it divides $f_0$. The derivative $f_0'=\sum_{e\in E\setminus\{e_0\}}(e-e_0)c_et^{e-e_0-1}$ has exactly $|E|-1$ nonzero coefficients, as $0<e-e_0<p$, and its exponents lie in $\{0,\dots,p-1\}$. If $(t-1)^{|E|}$ divided $f_0$, then $(t-1)^{|E|-1}$ would divide $f_0'$, contradicting the inductive hypothesis.

(b) We use the isomorphism $\FF_p\Omega_1\cong R=\FF_p[t]/(t-1)^p$ from the proof of Lemma~\ref{lem:uniserial}, under which $J^d(\Omega_1)$ corresponds to $(t-1)^dR$. An element of $R$ represented by a polynomial $f$ of degree less than $p$ lies in $(t-1)^dR$ if and only if $(t-1)^d$ divides $f$. Hence (a) shows that $J^{|X|}(\Omega_1)\cap\FF_pX=0$, and the assertion follows by comparing dimensions.
\end{proof}

\begin{remark}
In the language of coding theory, Lemma~\ref{lem:MDS}(b) says that the repeated-root cyclic codes $\langle(t-1)^d\rangle$ of length $p$ over $\FF_p$ are MDS codes; more generally, the minimum distances of all cyclic codes of length $p^s$ over finite fields of characteristic $p$ were determined by Dinh \cite{Dinh}. As its proof is very short, we have included it. Lemma~\ref{lem:MDS}(b) also says that for every subset $I$ of the coordinates of the base group $B$ of $C_p\wr C_p$, the subgroup $\prod_{i\in I}\langle x_i\rangle$ has a complement in $B$ that is normal in $C_p\wr C_p$.
\end{remark}

\begin{definition}\label{def:good}
Let $k\ge 0$. A subset $X\subseteq\Omega_k$ is \emph{good} if $\FF_p\Omega_k=J^{|X|}(\Omega_k)\oplus\FF_pX$, or equivalently, if $J^{|X|}(\Omega_k)\cap\FF_pX=0$.
\end{definition}

Clearly the property of being good is invariant under the action of $P_k$, and every subset of $\Omega_0$ is good. By Lemma~\ref{lem:MDS}(b), every subset of $\Omega_1$ is good. For $k\ge 2$ not every subset of $\Omega_k$ is good: if $X$ consists of two points of a minimal block, then $J^2(\Omega_k)\cap\FF_pX\ne 0$, as $\FF_pX$ contains a nonzero element of $\FF_p\Omega_k\cdot(g^{p^{k-1}}-1)\subseteq J^{p^{k-1}}(\Omega_k)$, where $g$ is a $p^k$-cycle in $P_k$.

\begin{lemma}\label{lem:goodcomplement}
Let $k\ge 0$ and $X\subseteq\Omega_k$. Then $X$ is good if and only if $\Omega_k\setminus X$ is good.
\end{lemma}

\begin{proof}
The standard bilinear form on $\FF_p\Omega_k$, for which $(e_\omega)$ is an orthonormal basis, is nondegenerate and $P_k$-invariant. Hence the orthogonal complement of $J^d$ is a submodule of dimension $d$, that is, $(J^d)^\perp=J^{p^k-d}$ by Lemma~\ref{lem:uniserial}. Moreover $(\FF_pX)^\perp=\FF_p(\Omega_k\setminus X)$. If $X$ is good, then $J^{|X|}+\FF_pX=\FF_p\Omega_k$, so $J^{p^k-|X|}\cap\FF_p(\Omega_k\setminus X)=0$; as $p^k-|X|=|\Omega_k\setminus X|$, this means that $\Omega_k\setminus X$ is good. The converse follows by symmetry.
\end{proof}

\begin{proposition}\label{prop:balancedgood}
Every balanced subset of $\Omega_k$ is good.
\end{proposition}

\begin{proof}
We argue by induction on $k$, the cases $k\le 1$ being settled above. Let $k\ge 2$ and let $X\subseteq\Omega_k$ be balanced. Fix a $p^k$-cycle $g\in P_k$ and $\omega_0\in\Omega_k$, and identify $\FF_p\Omega_k$ with $R=\FF_p[t]/(t-1)^{p^k}$ via $\omega_0g^x\mapsto t^x$ ($0\le x<p^k$). By Lemma~\ref{lem:uniserial}, $J^d(\Omega_k)$ corresponds to $u^dR$, where $u=t-1$. The top blocks are blocks for $\langle g\rangle$ of size $p^{k-1}$, so they are the orbits of $\langle g^p\rangle$; we number them so that $\Omega_k^{(j)}=\{\omega_0g^{j+py}:\ 0\le y<p^{k-1}\}$ for $0\le j<p$. The element $g^p$ stabilizes all top blocks, and its restriction to $\Omega_k^{(j)}$ is a $p^{k-1}$-cycle lying in the corresponding factor $E_j\cong P_{k-1}$ of the base group. Let $A=\FF_p[s]/(s-1)^{p^{k-1}}$ and $w=s-1\in A$, and identify $\FF_p\Omega_k^{(j)}$ with $A$ via $\omega_0g^{j+py}\mapsto s^y$. By Lemma~\ref{lem:uniserial}, applied to $E_j$, the submodule $J^e(\Omega_k^{(j)})$ corresponds to $w^eA$. Let $X_j\subseteq A$ be the set of monomials $s^y$ corresponding to $X\cap\Omega_k^{(j)}$; by induction, $\FF_pX_j\oplus w^{|X_j|}A=A$ for all $j$.

The assignment $s\mapsto t^p$ defines a ring homomorphism $A\to R$, as $(t^p-1)^{p^{k-1}}=(t-1)^{p^k}$; it maps $w$ to $u^p$. It is injective, because its kernel is an ideal of $A$, hence of the form $w^eA$, and $w^e\mapsto u^{pe}\ne 0$ for $e<p^{k-1}$. We regard $A$ as a subring of $R$. Since $u^p=w\in A$, the elements $1,u,\dots,u^{p-1}$ generate $R$ as an $A$-module, and as $\dim_{\FF_p}R=p\dim_{\FF_p}A$, they form a basis of $R$ over $A$. Finally, $\omega_0g^{j+py}$ corresponds to $t^js^y$.

Since $X$ is balanced, there is an integer $b\ge 0$ such that $|X_j|\in\{b,b+1\}$ for all $j$ and the set $I=\{j:\ |X_j|=b+1\}$ satisfies $\kappa:=|I|<p$. Then $a:=|X|=pb+\kappa$. As $u^au^i=u^{\kappa+i}w^b$, which equals $u^{\kappa+i-p}w^{b+1}$ if $\kappa+i\ge p$, we obtain
\[
u^aR=\bigoplus_{0\le i<\kappa}u^iw^{b+1}A\ \oplus\bigoplus_{\kappa\le i<p}u^iw^bA.
\]
Let $f\in\FF_pX\cap u^aR$, and write $f=\sum_{j=0}^{p-1}t^jf_j$ with $f_j\in\FF_pX_j\subseteq A$. Since $t^j=(1+u)^j=\sum_i\binom ji u^i$, we have $f=\sum_{i=0}^{p-1}u^i\alpha_i$ with $\alpha_i=\sum_j\binom ji f_j\in A$. By the decomposition of $u^aR$, we have $\alpha_i\in w^bA$ for all $i$ and $\alpha_i\in w^{b+1}A$ for $i<\kappa$. As the matrix $\bigl(\binom ji\bigr)_{0\le i,j<p}$ is unitriangular, the first condition yields $f_j\in w^bA$ for all $j$.

If $j\notin I$, then $f_j\in\FF_pX_j\cap w^{|X_j|}A=0$. If $j\in I$, then $\FF_pX_j$ maps isomorphically onto $A/w^{b+1}A$, so the elements of $\FF_pX_j\cap w^bA$ are the multiples of an element $\phi_j\in\FF_pX_j$ with $\phi_j\equiv w^b\pmod{w^{b+1}A}$; hence $f_j=\lambda_j\phi_j$ with $\lambda_j\in\FF_p$. For $i<\kappa$ we now obtain
\[
\alpha_i\equiv\Bigl(\sum_{j\in I}\binom ji\lambda_j\Bigr)w^b\pmod{w^{b+1}A},
\]
and therefore $\sum_{j\in I}\binom ji\lambda_j=0$ for $0\le i<\kappa$. This means that $(t-1)^\kappa$ divides the polynomial $h=\sum_{j\in I}\lambda_jt^j=\sum_iu^i\sum_{j\in I}\binom ji\lambda_j$. If $h\ne 0$, let $E\subseteq I$ be the set of exponents of $h$; then $|E|\le|I|=\kappa$, so $(t-1)^{|E|}$ divides $h$, contradicting Lemma~\ref{lem:MDS}(a). Hence $h=0$, so all $\lambda_j$ vanish and $f=0$. Thus $J^{|X|}(\Omega_k)\cap\FF_pX=0$, that is, $X$ is good.
\end{proof}

\subsection{The subgroups \texorpdfstring{$T_\Delta$}{T Delta}}\label{sec:TDelta}

Let $m\ge 2$. We define the \emph{block sum map}
\[
\rho\colon\FF_p\Omega_{m-1}\to\FF_p\Omega_{m-2},\qquad e_\omega\mapsto e_i\quad\text{for }\omega\in\beta_i.
\]

\begin{lemma}\label{lem:blocksum}
The map $\rho$ is a surjective homomorphism of $\FF_pP_{m-1}$-modules, and
\[
J^c(\Omega_{m-1})=\rho^{-1}\bigl(J^c(\Omega_{m-2})\bigr)\qquad\text{for }0\le c\le p^{m-2}.
\]
In particular, $H_{m-2}$ acts trivially on $H_{m-1}/J^c(\Omega_{m-1})$ for $0\le c\le p^{m-2}$.
\end{lemma}

\begin{proof}
The first assertion is clear, as $P_{m-1}$ permutes the minimal blocks. Hence $\rho^{-1}(J^c(\Omega_{m-2}))$ is a $P_{m-1}$-submodule of $\FF_p\Omega_{m-1}$ of codimension $c$, so it equals $J^c(\Omega_{m-1})$ by Lemma~\ref{lem:uniserial}. The last assertion follows because $H_{m-2}$ acts trivially on $\FF_p\Omega_{m-2}\cong\FF_p\Omega_{m-1}/\ker\rho$.
\end{proof}

For $\Delta\subseteq\Omega_{m-2}$ we put $c_\Delta=p^{m-2}-|\Delta|$ and define
\[
N_\Delta=J^{|\Delta|}(\Omega_{m-2})\le H_{m-2},\qquad M_\Delta=J^{c_\Delta}(\Omega_{m-1})\le H_{m-1},\qquad T_\Delta=Q\,N_\Delta\,M_\Delta,
\]
where $Q=H_0\cdots H_{m-3}$ as in Section~\ref{sec:orbits}. By Lemma~\ref{lem:uniserial}, $N_\Delta$ is normal in $P_{m-1}$ and $M_\Delta$ is normal in $P_m$; hence $T_\Delta$ is a subgroup of $P_m$. It depends only on $|\Delta|$, and its order is
\[
|T_\Delta|=p^{(p^{m-2}-1)/(p-1)}\cdot p^{p^{m-2}-|\Delta|}\cdot p^{p^{m-1}-c_\Delta}=p^{x_m}.
\]
Moreover, $T_\Delta\cap H_{m-1}=M_\Delta$, and the image of $T_\Delta$ in $P_{m-1}$ is $QN_\Delta$.

\begin{theorem}\label{thm:TDelta}
Let $m\ge 2$, let $\Delta\subseteq\Omega_{m-2}$ lie in a regular orbit of $P_{m-2}$ on $\Pow(\Omega_{m-2})$, and let $O_\Delta$ be the corresponding largest orbit of $P_m$ on $V_m$.
\begin{enumerate}[\rm(a)]
\item $O_\Delta$ is a regular orbit of $T_\Delta$ if and only if $\Delta$ is good.
\item If $\Delta$ is good, then the regular orbits of $T_\Delta$ are precisely the orbits $O_{\Delta'}$, where $\Delta'$ lies in a regular orbit of $P_{m-2}$, $\Delta'$ is good, and $|\Delta'|=|\Delta|$.
\item If $m\ge 3$, then $V_m$ is an irreducible $\FF_2T_\Delta$-module.
\end{enumerate}
\end{theorem}

\begin{proof}
(a) Let $L=H_{m-2}H_{m-1}$. Since $N_\Delta M_\Delta\le L$, $Q\cap L=1$ and $T_\Delta=N_\Delta M_\Delta Q$, Dedekind's law gives $T_\Delta\cap L=N_\Delta M_\Delta$. By Theorem~\ref{thm:largest}(d) we have $\C_{P_m}(v_\Delta)=Y_\Delta Z_\Delta\le L$ with $Y_\Delta\le H_{m-2}$ and $Z_\Delta\le H_{m-1}$. As every element of $L$ can be written uniquely as a product of an element of $H_{m-2}$ and an element of $H_{m-1}$, we obtain
\[
T_\Delta\cap\C_{P_m}(v_\Delta)=(N_\Delta\cap Y_\Delta)(M_\Delta\cap Z_\Delta).
\]
In additive notation, $Y_\Delta=\FF_p\Delta\le\FF_p\Omega_{m-2}$ and $Z_\Delta=\FF_p\Sigma_\Delta\le\FF_p\Omega_{m-1}$, where $\Sigma_\Delta=\{\sigma_i:\ i\notin\Delta\}$. Thus $N_\Delta\cap Y_\Delta=0$ if and only if $\Delta$ is good. Moreover, $\rho$ maps $\FF_p\Sigma_\Delta$ isomorphically onto $\FF_p(\Omega_{m-2}\setminus\Delta)$, and $M_\Delta=\rho^{-1}(J^{c_\Delta}(\Omega_{m-2}))$ by Lemma~\ref{lem:blocksum}. Hence $M_\Delta\cap Z_\Delta=0$ if and only if $J^{c_\Delta}(\Omega_{m-2})\cap\FF_p(\Omega_{m-2}\setminus\Delta)=0$, that is, if and only if $\Omega_{m-2}\setminus\Delta$ is good, which by Lemma~\ref{lem:goodcomplement} holds if and only if $\Delta$ is good. Now (a) follows from Lemma~\ref{lem:regular}.

(b) By Theorem~\ref{thm:largest} and Lemma~\ref{lem:regular}, the regular orbits of $T_\Delta$ are among the orbits $O_{\Delta'}$, where $\Delta'$ lies in a regular orbit of $P_{m-2}$, and as in (a), $O_{\Delta'}$ is a regular orbit of $T_\Delta$ if and only if $N_\Delta\cap\FF_p\Delta'=0$ and $M_\Delta\cap\FF_p\Sigma_{\Delta'}=0$. If $|\Delta'|>|\Delta|$, then $\dim N_\Delta+\dim\FF_p\Delta'>p^{m-2}$, so the first condition fails. If $|\Delta'|<|\Delta|$, then $\dim M_\Delta+\dim\FF_p\Sigma_{\Delta'}=p^{m-1}-c_\Delta+c_{\Delta'}>p^{m-1}$, so the second condition fails. If $|\Delta'|=|\Delta|$, then $T_{\Delta'}=T_\Delta$, and the assertion follows from (a), applied to $\Delta'$.

(c) Let $m\ge 3$. As $P_{m-2}\ne 1$, the set $\Omega_{m-2}$ does not lie in a regular orbit of $P_{m-2}$, so $|\Delta|<p^{m-2}$ and $N_\Delta\ne 0$. By Lemma~\ref{lem:uniserial}, $N_\Delta$ contains the element $\prod_iy_i$, which acts as a $p$-cycle on every minimal block. As $Q$ permutes the minimal blocks transitively, the image $QN_\Delta$ of $T_\Delta$ in $P_{m-1}$ acts transitively on $\Omega_{m-1}$. Moreover, $T_\Delta\cap H_{m-1}=M_\Delta$ contains $\ker\rho=J^{p^{m-2}}(\Omega_{m-1})$, which is the group described in Lemma~\ref{lem:irred}. Hence $V_m$ is irreducible for $T_\Delta$ by Lemma~\ref{lem:irred}.
\end{proof}

\begin{corollary}\label{cor:TUTD}
Let $m\ge 2$, and let $T_U=T_\Delta$ with $|\Delta|=(p^{m-2}+1)/2$ and $T_D=T_\Delta$ with $|\Delta|=(p^{m-2}-1)/2$. Explicitly,
\[
T_U=H_0\cdots H_{m-3}\,N_U\,M_U\qquad\text{and}\qquad T_D=H_0\cdots H_{m-3}\,N_D\,M_D,
\]
where $N_U,N_D\le H_{m-2}$ are normal in $H_0\cdots H_{m-2}$ with $|N_U|=p^{(p^{m-2}-1)/2}$ and $|N_D|=p^{(p^{m-2}+1)/2}$, and $M_U,M_D\le H_{m-1}$ are normal in $P_m$ with $|M_U|=p^{p^{m-1}-(p^{m-2}-1)/2}$ and $|M_D|=p^{p^{m-1}-(p^{m-2}+1)/2}$.
\begin{enumerate}[\rm(a)]
\item If $v\in V_m$ is of type $U_m$ (resp.\ $D_m$), then $v^{P_m}$ is a regular orbit of $T_U$ (resp.\ $T_D$).
\item If $m\ge 3$, then $V_m$ is an irreducible module for $T_U$ and for $T_D$. If $m=2$, then $T_U=H_1$ is the base group of $P_2$, which acts reducibly on $V_2$, and $T_D=M_0$ in the notation of Section~\ref{sec:m2}.
\item If $p=3$, then $O_U$ is the unique regular orbit of $T_U$, and $O_D$ is the unique regular orbit of $T_D$.
\end{enumerate}
\end{corollary}

\begin{proof}
(a) By Lemma~\ref{lem:types}, $v$ lies in a largest orbit, and $\Delta(v)$ is balanced of the appropriate size; so $\Delta(v)$ is good by Proposition~\ref{prop:balancedgood}, and (a) follows from Theorem~\ref{thm:TDelta}(a). Part (b) follows from Theorem~\ref{thm:TDelta}(c) and Section~\ref{sec:m2}. (c) If $p=3$, then by Lemma~\ref{lem:types} the largest orbits of $P_m$ are $O_U$ and $O_D$, and the corresponding sets $\Delta$ have different sizes. Thus (c) follows from (a) and Theorem~\ref{thm:TDelta}(b).
\end{proof}

\begin{remark}\label{rem:generalp}
Let $p$ be an arbitrary Mersenne prime and $m=3$. Then every subset of $\Omega_1$ is good, and the largest orbits of $P_3$ on $V_3$ correspond to the $P_1$-orbits of proper nonempty subsets of $\Omega_1$. By Theorem~\ref{thm:TDelta}(b), $T_\Delta$ has exactly $\binom{p}{|\Delta|}/p$ regular orbits. Hence $T_\Delta$ has a unique regular orbit if and only if $|\Delta|\in\{1,p-1\}$. In particular, for $p\ge 7$ the groups $T_U$ and $T_D$ have more than one regular orbit (for $p=7$ they have exactly five), whereas for $|\Delta|=1$ the group $T_\Delta$ is an irreducible subgroup of $P_3$ of order $p^{x_3}$ with a unique regular orbit. Together with Corollary~\ref{cor:m2} and the case $m=1$, this completes the proof of Theorem~\ref{thmC}.
\end{remark}

\subsection{Passing from \texorpdfstring{$P_{m+1}$}{P(m+1)} to \texorpdfstring{$P_m$}{P(m)}}\label{sec:descent}

The following result shows that the existence of subgroups with a unique regular orbit of the largest possible size propagates from $P_{m+1}$ down to $P_m$.

\begin{proposition}\label{prop:descent}
Let $m\ge 2$. If some subgroup of $P_{m+1}$ has a largest orbit of $P_{m+1}$ on $V_{m+1}$ as its unique regular orbit, then some subgroup of $P_m$ has a largest orbit of $P_m$ on $V_m$ as its unique regular orbit.
\end{proposition}

\begin{proof}
Let $V_{m+1}=X_1\oplus\dots\oplus X_p$ and $E=E_1\times\dots\times E_p$ be as in Definition~\ref{def:types}, so that $X_j$ is the natural module for $E_j\cong P_m$. Let $v=v_1+\dots+v_p\in V_{m+1}$ with $v_j\in X_j$. Then $\Delta(v)$ is the disjoint union of the sets $\Delta(v_j)$. Hence, by Theorem~\ref{thm:largest}, applied to $P_{m+1}$ and to $P_m$, and by Lemma~\ref{lem:powerset}(a), the vector $v$ lies in a largest orbit of $P_{m+1}$ if and only if each $v_j$ lies in a largest orbit of $E_j$ and these orbits, identified via the top layer of $P_{m+1}$, are not all equal. Moreover, two such vectors lie in the same $P_{m+1}$-orbit if and only if their tuples of $E_j$-orbits agree up to a cyclic permutation, and in this case $\C_{P_{m+1}}(v)\le H_{m-1}H_m\le E$.

Let $T\le P_{m+1}$ have the largest orbit $O$ as its unique regular orbit, so that $|T|=|O|=p^{x_{m+1}}$, and let $v=v_1+\dots+v_p\in O$. Since $T$ is transitive on $O$ and $\C_{P_{m+1}}(v)\le E$, we have $TE=P_{m+1}$. Hence $|T\cap E|=|T|/p=|v^E|$, and $T\cap E$ acts regularly on $v^E=v_1^{E_1}\times\dots\times v_p^{E_p}$. Fix $j$ and let $T_j$ be the set of elements of $T\cap E$ that fix $v_i$ for all $i\ne j$. Then $T_j$ acts regularly on the set of vectors $w+\sum_{i\ne j}v_i$ with $w\in v_j^{E_j}$. Moreover, $T_j$ acts faithfully on $X_j$, since an element of $T_j$ acting trivially on $X_j$ fixes $v$. Thus $T_j$ induces a subgroup $S_j\le E_j\cong P_m$ of order $|v_j^{E_j}|=p^{x_m}$ having $v_j^{E_j}$ as a regular orbit.

Suppose that $S_j$ has another regular orbit $O_j'$, which is a largest orbit of $E_j$ by Lemma~\ref{lem:regular}, and that the orbits $v_i^{E_i}$ with $i\ne j$ are not all equal to $O_j'$. Let $w\in O_j'$ and $v'=w+\sum_{i\ne j}v_i$. Then $v'$ lies in a largest orbit of $P_{m+1}$, which differs from $O$ because the multisets of the corresponding $E_i$-orbits differ. Every element of $T\cap\C_{P_{m+1}}(v')$ lies in $T_j$ and fixes $w$, hence is trivial. Thus $v'^{P_{m+1}}$ is a second regular orbit of $T$ by Lemma~\ref{lem:regular}, a contradiction.

Finally, the tuple $(v_1^{E_1},\dots,v_p^{E_p})$ is not constant. If it consists of one orbit occurring once and another orbit occurring $p-1\ge 2$ times, choose $j$ such that $v_j^{E_j}$ is the latter orbit; otherwise let $j$ be arbitrary. In both cases the orbits $v_i^{E_i}$ with $i\ne j$ are not all equal, so the largest orbit $v_j^{E_j}$ of $E_j\cong P_m$ is the unique regular orbit of $S_j$.
\end{proof}

\begin{corollary}\label{cor:descent}
Let $p$ be a Mersenne prime. The set of integers $m\ge 1$ for which some subgroup of $P_m$ has a largest orbit of $P_m$ on $V_m$ as its unique regular orbit is an initial segment of the positive integers. It contains $\{1,2,3\}$, and it contains all positive integers if $p=3$. In particular, if for some $m_0\ge 3$ no subgroup of $P_{m_0}$ has a largest orbit of $P_{m_0}$ as its unique regular orbit, then the same holds for all $m\ge m_0$.
\end{corollary}

\begin{proof}
By Proposition~\ref{prop:descent}, the set in question is an initial segment, and the last assertion is a reformulation of this fact. The set contains $1$, $2$ and $3$ by the case $m=1$, Corollary~\ref{cor:m2} and Remark~\ref{rem:generalp}, respectively, and for $p=3$ it contains every $m\ge 2$ by Corollary~\ref{cor:TUTD}(c).
\end{proof}

\subsection{Non-uniqueness}\label{sec:nonunique}

Recall the notation $\mathcal{S}(O)$ introduced before Lemma~\ref{lem:image}. For $p=3$ and $m\ge 2$ we fix subsets $\Delta_U,\Delta_D\subseteq\Omega_{m-2}$ with $O_U=O_{\Delta_U}$ and $O_D=O_{\Delta_D}$, so that $|\Delta_U|=(3^{m-2}+1)/2$, $|\Delta_D|=(3^{m-2}-1)/2$, $T_U=T_{\Delta_U}$ and $T_D=T_{\Delta_D}$ (Lemma~\ref{lem:types} and Corollary~\ref{cor:TUTD}); moreover, $N_U=N_{\Delta_U}$ and $M_U=M_{\Delta_U}$.

\begin{proposition}\label{prop:nonunique}
Let $p=3$ and $m\ge 3$.
\begin{enumerate}[\rm(a)]
\item Let $T\le P_m$ with $|T|=3^{x_m}$. Assume that the image $\bar T$ of $T$ in $P_{m-1}$ acts regularly on $\mathcal{S}(O_U)$ and that $(T\cap H_{m-1})\cap\C_{H_{m-1}}(v)=1$ for some $v\in O_U$. Then $O_U$ is the unique regular orbit of $T$.
\item The hypotheses of \textup{(a)} are satisfied by every subgroup $T\le P_m$ with $TH_{m-1}=T_UH_{m-1}$ and $T\cap H_{m-1}=M_U$, and every such $T$ acts irreducibly on $V_m$.
\item There is a subgroup $T$ as in \textup{(b)} that is not conjugate to $T_U$ in $P_m$.
\end{enumerate}
\end{proposition}

\begin{proof}
(a) As $\bar T$ acts regularly on $\mathcal{S}(O_U)$, the image of $T\cap\C_{P_m}(v)$ in $P_{m-1}$ is contained in $\bar T\cap\Stab_{P_{m-1}}(\supp(v))=1$, so $T\cap\C_{P_m}(v)\le(T\cap H_{m-1})\cap\C_{H_{m-1}}(v)=1$. By Lemma~\ref{lem:regular}, $O_U$ is a regular orbit of $T$. If $O_D$ were a regular orbit of $T$, then $\bar T$ would act transitively on $\mathcal{S}(O_D)$, which has $|P_{m-1}|/3^{|\Delta_D|}$ elements. But $|\bar T|=|\mathcal{S}(O_U)|=|P_{m-1}|/3^{|\Delta_U|}$ and $|\Delta_U|=|\Delta_D|+1$, a contradiction.

(b) Let $T$ be as in (b), so that $|T|=|QN_U||M_U|=|T_U|=3^{x_m}$, and put $\Delta=\Delta_U$. Since $T_U$ acts regularly on $O_U$ by Corollary~\ref{cor:TUTD}, its image $QN_\Delta$ in $P_{m-1}$ acts transitively on $\mathcal{S}(O_U)$, and as $|QN_\Delta|=|P_{m-1}|/3^{|\Delta|}=|\mathcal{S}(O_U)|$, it acts regularly. Moreover, $M_\Delta\cap\C_{H_{m-1}}(v_\Delta)=M_\Delta\cap Z_\Delta=1$ by the proof of Theorem~\ref{thm:TDelta}(a). The irreducibility follows from Lemma~\ref{lem:irred} exactly as in the proof of Theorem~\ref{thm:TDelta}(c).

(c) Put $\bar T_U=QN_U$ and $X=H_{m-1}/M_U$, a $\bar T_U$-module on which $N_U$ acts trivially by Lemma~\ref{lem:blocksum}. The subgroups $T$ as in (b) correspond to the complements of $X$ in $H_{m-1}\bar T_U/M_U=X\rtimes\bar T_U$, that is, to the cocycles $\delta\in Z^1(\bar T_U,X)$ via $T/M_U=\{t\delta(t):\ t\in\bar T_U\}$; and $T_U$ corresponds to $\delta=0$. As $m\ge 3$, we have $N_U\ne 0$ and $c:=3^{m-2}-|\Delta_U|\ge 1$, so $X_0=J^{c-1}(\Omega_{m-1})/M_U$ is a one-dimensional subgroup of $X$ on which $P_m$ acts trivially. Let $\psi\colon N_U\to X_0$ be the composition of the natural map $N_U\to N_U/[N_U,Q]\cong C_3$ with an isomorphism $C_3\to X_0$. Here $[N_U,Q]=[N_U,P_{m-1}]$ has index $3$ in $N_U$ by Lemma~\ref{lem:uniserial}, and $\psi(n^q)=\psi(n)$ for all $n\in N_U$ and $q\in Q$. Since $\bar T_U=N_UQ$ with $N_U\cap Q=1$ and $N_U\trianglelefteq\bar T_U$, the map $\delta\colon\bar T_U\to X_0$, $nq\mapsto\psi(n)$ ($n\in N_U$, $q\in Q$), is a well-defined homomorphism; as $\bar T_U$ acts trivially on $X_0$, it is a cocycle with values in $X$, and its restriction to $N_U$ is $\psi\ne 0$. Let $T_\delta$ be the corresponding subgroup.

Suppose that $T_\delta^g=T_U$ for some $g\in P_m$. Then $g$ normalizes $T_UH_{m-1}=H_{m-1}\bar T_U$, so $g=\tau d$ with $\tau\in\N_{P_{m-1}}(\bar T_U)$ and $d\in H_{m-1}$. Conjugation by $\tau$ maps $T_\delta$ to the subgroup corresponding to the cocycle $t\mapsto\delta(t^{\tau^{-1}})^{\tau}$, whose restriction to $N_U$ is again nonzero, and conjugation by $d$ changes a cocycle by a coboundary, whose restriction to $N_U$ vanishes as $N_U$ acts trivially on $X$. Hence $T_\delta^g$ corresponds to a cocycle with nonzero restriction to $N_U$, so $T_\delta^g\ne T_U$, a contradiction.
\end{proof}

Proposition~\ref{prop:nonunique} proves Theorem~\ref{thmD}(d). In fact, more is true for small $m$: a computation with GAP shows that for $m=3$ and $m=4$ the group $T_\delta$ constructed in the proof of (c) is not even isomorphic to $T_U$, as $T_U/T_U'$ is elementary abelian of rank $4$ (resp.\ $5$), whereas $T_\delta/T_\delta'$ is elementary abelian of rank $3$ (resp.\ $4$).

\subsection{The case \texorpdfstring{$p=3$}{p=3} and \texorpdfstring{$m=3$}{m=3}}\label{sec:m3}

In this subsection let $p=3$, so that $F=\FF_4$, $|P_3|=3^{13}$, $|V_3|=4^9$ and $x_3=10$. Recall that $P_2=H_0H_1$ and that $M_0,M_1,M_2$ denote the maximal subgroups of $P_2$ different from $H_1$ (Lemma~\ref{lem:Pstructure}).

\begin{theorem}\label{thm:m3}
Let $p=3$. Up to conjugacy in $P_3$, there are exactly $30$ subgroups $T\le P_3$ of order $3^{10}$ having exactly one regular orbit on $V_3$; for $15$ of them this orbit is $O_U$, and for the other $15$ it is $O_D$. All of them act irreducibly on $V_3$. In each case $T\cap H_2=J^c(\Omega_2)$ with $c\in\{1,2,3\}$, and under $\GL(V_3)$-conjugacy the $30$ groups fall into the $13$ classes listed in Table~\ref{tab:m3}. In particular, there are exactly $13$ isomorphism types of pairs $(G,V)$, where $G$ is a $3$-group of order $3^{10}$ and $V$ is a faithful irreducible $G$-module of order $4^9$ such that $G$ has exactly one regular orbit on $V$.
\end{theorem}

\begin{table}[ht]
\centering
\footnotesize
\setlength{\tabcolsep}{5pt}
\begin{tabular}{@{}rlcccccrr@{}}
\toprule
No. & $\bar T$ & $|H_2:M|$ & $T=\bar TM$ & regular orbit & $P_3$-classes & $T/T'$ & $\exp T$ & orbits\\
\midrule
1 & $C_3\times C_3$ & $3$ & yes ($T_U$) & $O_U$ & $1$ & $C_3^4$ & $9$ & $66$\\
2 & $C_3\times C_3$ & $3$ & no & $O_U$ & $4$ & $C_3^3$ & $9$ & $64$\\
3 & $C_9$ & $3$ & yes & $O_U$ & $2$ & $C_3\times C_9$ & $9$ & $62$\\
4 & $C_9$ & $3$ & no & $O_U$ & $4$ & $C_3\times C_9$ & $27$ & $60$\\
5 & $M_0$ & $3^2$ & yes ($T_D$) & $O_D$ & $1$ & $C_3^4$ & $9$ & $42$\\
6 & $M_0$ & $3^2$ & no & $O_D$ & $2$ & $C_3\times C_9$ & $9$ & $40$\\
7 & $M_0$ & $3^2$ & no & $O_D$ & $2$ & $C_3^2$ & $9$ & $36$\\
8 & $M_1$, $M_2$ & $3^2$ & yes & $O_D$ & $2$ & $C_3^3$ & $9$ & $38$\\
9 & $M_1$, $M_2$ & $3^2$ & no & $O_D$ & $4$ & $C_3\times C_9$ & $27$ & $36$\\
10 & $P_2$ & $3^3$ & no & $O_D$ & $2$ & $C_3^3$ & $9$ & $34$\\
11 & $P_2$ & $3^3$ & no & $O_U$ & $2$ & $C_3^2$ & $9$ & $34$\\
12 & $P_2$ & $3^3$ & no & $O_D$ & $2$ & $C_3^3$ & $27$ & $28$\\
13 & $P_2$ & $3^3$ & no & $O_U$ & $2$ & $C_3^3$ & $27$ & $32$\\
\bottomrule
\end{tabular}
\medskip
\caption{The subgroups $T\le P_3$ of order $3^{10}$ with exactly one regular orbit on $V_3$ ($p=3$), up to $\GL(V_3)$-conjugacy. Here $M=T\cap H_2=J^c(\Omega_2)$ with $3^c=|H_2:M|$, and $\bar T=TH_2/H_2\le P_2$ is given up to conjugacy in $P_2$; the subgroups $C_3\times C_3$ and $C_9$ are regular on $\Omega_2$. The fourth column indicates whether $T$ is the semidirect product of $M$ and the subgroup $\bar T$ of $P_2=H_0H_1$, the sixth column gives the number of $P_3$-conjugacy classes contained in the $\GL(V_3)$-class, and the last column gives the number of orbits of $T$ on $V_3$.}
\label{tab:m3}
\end{table}

\begin{proof}[Proof (computer-assisted)]
Let $T\le P_3$ with $|T|=3^{10}$, and let $\bar T=TH_2/H_2$ and $M=T\cap H_2$, so that $|\bar T||M|=3^{10}$. If $T$ has a regular orbit $O\in\{O_U,O_D\}$, then $T$ acts transitively on $O$, and so $\bar T$ acts transitively on $\mathcal{S}(O)$, which has $9$ elements if $O=O_U$ and $27$ elements if $O=O_D$. Among the $20$ conjugacy classes of subgroups of $P_2$, exactly $7$ have this property for $O_U$ or for $O_D$: the three classes of regular subgroups of order $9$ (one isomorphic to $C_3\times C_3$ and two cyclic ones), the classes of $M_0$, $M_1$, $M_2$, and $P_2$ itself. In accordance with Lemma~\ref{lem:image}, these are exactly the subgroups $\bar T\le P_2$ (up to conjugacy) with $\bar TH_1=P_2$ and $\bar T\cap H_1=J^c(\Omega_1)$ for some $c\le 2$. For each such $\bar T$ we computed, using the MeatAxe, all $\bar T$-submodules $M$ of $H_2\cong\FF_3\Omega_2$ of dimension $10-\log_3|\bar T|$ (in each case there is exactly one, namely $J^c(\Omega_2)$ with $c=\log_3|\bar T|-1$), and then all complements of $H_2/M$ in $H_2\bar T/M$ up to conjugacy. For each resulting group $T$ we decided whether $T\cap\C_{P_3}(v)=1$ for $v=v_U$ and for $v=v_D$, where $v_U\in O_U$ and $v_D\in O_D$ are as in Theorem~\ref{thm:largest}(d); by Lemma~\ref{lem:regular}, this determines the regular orbits of $T$. The resulting $30$ classes of groups with exactly one regular orbit were then checked independently by computing all orbits of each of them on $V_3$. Irreducibility was verified with the MeatAxe; it also follows from Lemma~\ref{lem:irred}, as each $\bar T$ is transitive on $\Omega_2$ and $T\cap H_2\supseteq J^3(\Omega_2)=\ker\rho$.

To determine the $\GL(V_3)$-classes, let $\mathcal{W}=\GammaL_1(4)\wr S_9\le\GL(V_3)$ be the group of $\FF_2$-linear monomial maps permuting the subspaces $Fe_\omega$, $\omega\in\Omega_2$; viewing $P_3$ as a Sylow $3$-subgroup of $S_{27}$ acting on the set $\{\lambda e_\omega:\ \lambda\in F^\times,\ \omega\in\Omega_2\}$, the group $\mathcal{W}$ is the stabilizer in $S_{27}$ of the system of minimal blocks. Conjugacy in $\mathcal{W}$ fuses the $30$ classes into the $13$ classes of Table~\ref{tab:m3}. Since the number of orbits on $V_3$ and the isomorphism type of $T/T'$ are invariant under conjugacy in $\GL(V_3)$, and these invariants separate the $13$ classes, there is no further fusion in $\GL(V_3)$. The last assertion follows from Theorem~\ref{thm:reduction} and Remark~\ref{rem:sylow}.
\end{proof}

\begin{remark}
\begin{enumerate}[(a)]
\item The groups in rows~1--4 of Table~\ref{tab:m3} are exactly the groups described in Proposition~\ref{prop:nonunique}(a) for $m=3$: here $\bar T$ is regular on $\Omega_2$, $T\cap H_2=J^1(\Omega_2)$ is the augmentation submodule, and every complement of $H_2/(T\cap H_2)$ in $H_2\bar T/(T\cap H_2)$ gives a group with exactly one regular orbit.
\item In rows~5--9 the group $\bar T$ acts regularly on $\mathcal{S}(O_D)$, so all groups obtained from these $\bar T$ as in the proof have $O_D$ as a regular orbit; among the $45$ complement classes arising in the computation, exactly $15$ give groups that do not also have $O_U$ as a regular orbit.
\item In rows~10--13 we have $\bar T=P_2$, which does not act regularly on $\mathcal{S}(O_U)$ or $\mathcal{S}(O_D)$. Thus the groups with exactly one regular orbit do not all arise as in Proposition~\ref{prop:nonunique}(a) or its analogue for $O_D$.
\item The computation also shows that there are $36$ conjugacy classes of subgroups of $P_3$ of order $3^{10}$ having two regular orbits on $V_3$.
\end{enumerate}
\end{remark}

Corollary~\ref{cor:m2} and Theorem~\ref{thm:m3} prove parts (a) and (b) of Theorem~\ref{thmD}.

\subsection{The case \texorpdfstring{$p=3$}{p=3} and \texorpdfstring{$m=4$}{m=4}}\label{sec:m4}

Now let $p=3$ and $m=4$, so that $|P_4|=3^{40}$, $|V_4|=4^{27}$ and $x_4=31$. Lemma~\ref{lem:image} makes the method of Theorem~\ref{thm:m3} applicable in this case.

\begin{theorem}\label{thm:m4}
Let $p=3$. Up to conjugacy in $P_4$, there are exactly $4914$ subgroups $T\le P_4$ of order $3^{31}$ having exactly one regular orbit on $V_4$; for $2457$ of them this orbit is $O_U$, and for the other $2457$ it is $O_D$. All of them act irreducibly on $V_4$. For each of them there is $c\in\{0,\dots,5\}$ such that $\bar T\cap H_2=J^c(\Omega_2)$ and $T\cap H_3=J^{9-c}(\Omega_3)$, where $\bar T=TH_3/H_3$, and the numbers of conjugacy classes for the various values of $c$ are as given in Table~\ref{tab:m4}.
\end{theorem}

\begin{table}[ht]
\centering
\footnotesize
\begin{tabular}{@{}lrrrrrrr@{}}
\toprule
$c$ & $0$ & $1$ & $2$ & $3$ & $4$ & $5$ & total\\
\midrule
classes of images $\bar T$ & $1$ & $9$ & $27$ & $27$ & $27$ & $27$ & $118$\\
classes of $T$ with regular orbit $O_U$ & $162$ & $360$ & $54$ & $396$ & $0$ & $1485$ & $2457$\\
classes of $T$ with regular orbit $O_D$ & $162$ & $162$ & $648$ & $396$ & $1089$ & $0$ & $2457$\\
\bottomrule
\end{tabular}
\medskip
\caption{The conjugacy classes of subgroups $T\le P_4$ of order $3^{31}$ with exactly one regular orbit on $V_4$ ($p=3$), according to the integer $c$ with $\bar T\cap H_2=J^c(\Omega_2)$.}
\label{tab:m4}
\end{table}

\begin{proof}[Proof (computer-assisted)]
Let $T\le P_4$ have order $3^{31}$ and a regular orbit $O_\Delta$. By Lemma~\ref{lem:image}, the image $\bar T$ of $T$ in $P_3$ satisfies $\bar TH_2=P_3$ and $\bar T\cap H_2=J^c(\Omega_2)$ with $c\le|\Delta|\le 5$. Hence $\bar T$ is the preimage of a complement of $H_2/J^c(\Omega_2)$ in $P_3/J^c(\Omega_2)$; up to conjugacy in $P_3$ there are $118$ such groups, distributed over $c$ as in Table~\ref{tab:m4}. For each of them, the MeatAxe shows that $J^{9-c}(\Omega_3)$ is the only $\bar T$-submodule of $H_3$ of dimension $31-\log_3|\bar T|=18+c$, so that $M:=T\cap H_3=J^{9-c}(\Omega_3)$ (this also follows from Proposition~\ref{prop:Mforced} below). The groups $T$ with given $\bar T$ and $M$ correspond to the complements of $H_3/M$ in $H_3\bar T/M$, and their conjugacy classes under $H_3\bar T$ correspond to the elements of $H^1(\bar T,H_3/M)$, which we computed with the GAP function \texttt{OneCocycles}. For each of these classes we decided, as in the proof of Theorem~\ref{thm:m3}, whether $T\cap\C_{P_4}(v)=1$ for $v=v_U$ and for $v=v_D$. If $T$ and $T'$ are such groups and $T^g=T'$ with $g\in P_4$, then their images in $P_3$ are conjugate by the image of $g$; and if $T$ and $T'$ have the same image $\bar T$, then $g\in H_3\N_{P_3}(\bar T)$. Thus the conjugacy classes of the groups $T$ with image $\bar T$ correspond to the orbits of $\N_{P_3}(\bar T)$ on the relevant elements of $H^1(\bar T,H_3/M)$; here $|\N_{P_3}(\bar T):\bar T|=3$ unless $\bar T=P_3$. Finally, $\bar T$ acts transitively on $\Omega_3$, since it maps onto $P_2$ and contains $J^c(\Omega_2)\ni\prod_iy_i$, and $T\cap H_3=J^{9-c}(\Omega_3)$ contains $\ker\rho=J^9(\Omega_3)$; so $T$ acts irreducibly on $V_4$ by Lemma~\ref{lem:irred}. As a check, the same program for $m=3$ reproduces the $30$ classes of Theorem~\ref{thm:m3}, which were obtained there by testing conjugacy directly.
\end{proof}

The computation also shows that $11934$ elements of the cohomology groups considered in the proof give subgroups with two regular orbits. We have not determined how the $4914$ classes of Theorem~\ref{thm:m4} fuse under conjugation in $\GL(V_4)$: the method used for $m=3$ is not available, as $|V_4|=4^{27}$ is too large for computing all orbits on $V_4$.

Theorem~\ref{thm:m4} proves part (c) of Theorem~\ref{thmD}.

\section{Mersenne primes \texorpdfstring{$p\ge 7$}{p>=7}}\label{sec:pge7}

In this section we prove Theorem~\ref{thmE}. By Corollary~\ref{cor:descent}, it suffices to show that no subgroup of $P_4$ has a largest orbit of $P_4$ on $V_4$ as its unique regular orbit. The proof only uses linear algebra over $\FF_p$, and it proceeds as follows. Let $T\le P_4$ be a subgroup of order $p^{x_4}$ with a regular orbit. In Sections~\ref{sec:coords}--\ref{sec:criterion} we show that $T\cap H_2H_3$ is essentially the graph of a linear map, and that the regular orbits of $T$ are detected by the invertibility of certain square matrices built from this map. In Section~\ref{sec:forms} we show that these matrices are circulant inside each block of $p$ points and constant between different blocks. As a consequence, changing the subset in a single block affects the invertibility in a very simple way, governed by a single ``forbidden value'' (Section~\ref{sec:oneblock}). Section~\ref{sec:circulants} contains a computation with small circulant matrices, and Section~\ref{sec:proofE} completes the proof. The assumption $p\ge 7$ is only used in Sections~\ref{sec:circulants} and~\ref{sec:proofE}, and there only in the form $p\ge 5$.

\subsection{Coordinates and types}\label{sec:coords}

We fix coordinates on $\Omega_1$ and $\Omega_2$ as follows. Let $\Omega_1=\ZZ_p$, where $P_1$ is generated by $b\mapsto b+1$, and let $\Omega_2=\ZZ_p\times\ZZ_p$, where for $b\in\ZZ_p$ the set $B_b=\{b\}\times\ZZ_p$ is the top block (and minimal block) of $\Omega_2$ corresponding to $b\in\Omega_1$. We may assume that $H_1=\langle h_b:\ b\in\ZZ_p\rangle$ and $H_0=\langle\hat h\rangle$, where $h_b$ maps $(b,x)$ to $(b,x+1)$ and fixes the points outside $B_b$, and $\hat h$ maps $(b,x)$ to $(b+1,x)$. We regard the elements of $\ZZ_p$ also as elements of $\FF_p$.

For $u\in\FF_p\Omega_2$ we write $u_{(b,x)}$ for the coefficient of $e_{(b,x)}$ in $u$, and we put
\[
\varepsilon_b(u)=\sum_{x\in\ZZ_p}u_{(b,x)}\qquad\text{and}\qquad \ell_b(u)=\sum_{x\in\ZZ_p}x\,u_{(b,x)}\qquad(b\in\ZZ_p).
\]
Note that $\ell_b(uh_b)=\ell_b(u)+\varepsilon_b(u)$. For $m=3$, the block sum map of Section~\ref{sec:TDelta} is $\rho_2\colon\FF_p\Omega_2\to\FF_p\Omega_1$, $u\mapsto\sum_b\varepsilon_b(u)e_b$. Its kernel is $K=\bigoplus_bK_b$, where $K_b$ denotes the set of vectors $u\in\FF_pB_b$ with $\varepsilon_b(u)=0$, and by Lemma~\ref{lem:blocksum} we have
\begin{equation}\label{eq:Jc}
J^c(\Omega_2)=\rho_2^{-1}\bigl(J^c(\Omega_1)\bigr)\supseteq K=J^p(\Omega_2)\qquad\text{for }0\le c\le p.
\end{equation}
In the notation of the proof of Lemma~\ref{lem:MDS}, a polynomial $f$ of degree less than $p$ is divisible by $(t-1)^2$ if and only if $f(1)=f'(1)=0$; hence $J^2(\Omega_1)$ consists of the vectors $v$ with $\sum_bv_b=\sum_bb\,v_b=0$. As $K_b\cong J^1(\Omega_1)$ as modules for $\langle h_b\rangle\cong P_1$, Lemma~\ref{lem:uniserial} yields
\begin{equation}\label{eq:Kb}
\{uh_b-u:\ u\in K_b\}=\{u\in K_b:\ \ell_b(u)=0\}.
\end{equation}

Let $\langle\cdot,\cdot\rangle$ be the standard bilinear form on $\FF_p\Omega_k$, as in the proof of Lemma~\ref{lem:goodcomplement}, so that $J^d(\Omega_k)^\perp=J^{p^k-d}(\Omega_k)$, and let $u\ast w=\sum_\omega u_\omega w_\omega e_\omega$ be the coordinatewise product; thus $\langle u\ast s,w\rangle=\langle u,s\ast w\rangle$.

A \emph{type} is an orbit of $P_1$ on the set of proper nonempty subsets of $\ZZ_p$; as $p$ is a prime, these are exactly the regular orbits of $P_1$ on $\Pow(\Omega_1)$. For each type $\tau$ we fix a representative $\Gamma_\tau\subseteq\ZZ_p$, and for a $p$-tuple $\mathbf t=(t_b)_{b\in\ZZ_p}$ of types we put
\[
\Delta(\mathbf t)=\bigcup_{b\in\ZZ_p}\{b\}\times\Gamma_{t_b}\subseteq\Omega_2.
\]
By Lemma~\ref{lem:powerset}(a) and its proof, every regular orbit of $P_2$ on $\Pow(\Omega_2)$ contains $\Delta(\mathbf t)$ for some non-constant tuple $\mathbf t$, and $\Delta(\mathbf t)$ and $\Delta(\mathbf t')$ lie in the same $P_2$-orbit if and only if $\mathbf t'$ is a cyclic shift of $\mathbf t$. Hence, by Theorem~\ref{thm:largest}, the largest orbits of $P_4$ on $V_4$ are the orbits $O_{\Delta(\mathbf t)}$ with $\mathbf t$ non-constant, and $O_{\Delta(\mathbf t)}=O_{\Delta(\mathbf t')}$ if and only if $\mathbf t'$ is a cyclic shift of $\mathbf t$. As there is only one type of $1$-element sets and only one type of $(p-1)$-element sets, we have $p+1\le|\Delta(\mathbf t)|\le p^2-p-1$ for every non-constant tuple $\mathbf t$.

\subsection{The intersection with \texorpdfstring{$H_3$}{H3}}\label{sec:Mforced}

\begin{lemma}\label{lem:products}
Let $0\le c<p^2$ and $0\ne s\in\FF_p\Omega_2$. Then $J^c(\Omega_2)\ast s\not\subseteq J^{c+1}(\Omega_2)$.
\end{lemma}

\begin{proof}
Fix a $p^2$-cycle $g\in P_2$ and $\omega_0\in\Omega_2$. For $0\le y<p^2$ put $e(y)=e_{\omega_0g^y}$ and write $y=y_0+py_1$ with $0\le y_0,y_1<p$. As in the proof of Lemma~\ref{lem:uniserial}, the map $e(y)\mapsto t^y$ identifies $\FF_p\Omega_2$ with $\FF_p[t]/(t-1)^{p^2}$ and $J^d(\Omega_2)$ with the ideal generated by $(t-1)^d$. Since $t^y=\sum_k\binom yk(t-1)^k$, a vector $\sum_ya_ye(y)$ lies in $J^d(\Omega_2)$ if and only if $\sum_y\binom yka_y=0$ for $0\le k<d$. Thus $J^d(\Omega_2)=\mathcal B_d^\perp$, where $\mathcal B_d$ is spanned by the vectors $\chi_k=\sum_y\binom yk e(y)$ with $0\le k<d$, and $\mathcal B_d=J^d(\Omega_2)^\perp=J^{p^2-d}(\Omega_2)$.

For $0\le i,j<p$ let $\mu_{i,j}=\sum_yy_0^iy_1^je(y)$ (with $0^0=1$), and call $i+pj$ the \emph{index} of $(i,j)$. By Lucas' theorem, $\binom yk\equiv\binom{y_0}{k_0}\binom{y_1}{k_1}\pmod p$ for $k=k_0+pk_1$ with $0\le k_0,k_1<p$, and $\binom{y_0}{k_0}$ is a polynomial in $y_0$ of degree $k_0$ with leading coefficient $1/k_0!$. As the set of pairs of index less than $d$ is closed under decreasing either entry, it follows that $\mathcal B_d$ is spanned by the $\mu_{i,j}$ with $i+pj<d$. The $p^2$ vectors $\mu_{i,j}$ form a basis of $\FF_p\Omega_2$, as the monomials $X_0^iX_1^j$ with $i,j<p$ form a basis of the space of maps $\FF_p^2\to\FF_p$. Moreover, $\mu_{i,j}\ast\mu_{i',j'}=\mu_{i\oplus i',\,j\oplus j'}$, where $i\oplus i'=i+i'$ if $i+i'<p$ and $i\oplus i'=i+i'-(p-1)$ otherwise, since $Y^p=Y$ for $Y\in\FF_p$.

Suppose that $J^c(\Omega_2)\ast s\subseteq J^{c+1}(\Omega_2)$. As $\langle u\ast s,w\rangle=\langle u,s\ast w\rangle$, this means that $s\ast\mathcal B_{c+1}\subseteq J^c(\Omega_2)^\perp=\mathcal B_c$. Since $\mu_{0,0}\in\mathcal B_{c+1}$, we obtain $s\in\mathcal B_c$. 
If $c=0$, this already contradicts $s\ne 0$. Hence we may assume $c\ge 1$. 
Write $s=\sum_{i,j}\lambda_{i,j}\mu_{i,j}$, let $(i_0,j_0)$ be the pair of largest index with $\lambda_{i_0,j_0}\ne 0$, so that $i_0+pj_0<c$, and write $c=c_0+pc_1$ with $0\le c_0,c_1<p$. Then $j_0\le c_1$; put $j'=c_1-j_0$, $i'=\max(c_0-i_0,0)$ and $c_0'=\max(c_0,i_0)$. The index of $(i',j')$ is at most $c$: if $j_0=c_1$ it is $c_0-i_0$, if $1\le j_0<c_1$ it is at most $p-1+p(c_1-1)<c$, and if $j_0=0$ it is at most $c_0+pc_1$. Hence $\mu_{i',j'}\in\mathcal B_{c+1}$. In $s\ast\mu_{i',j'}=\sum\lambda_{i,j}\mu_{i\oplus i',\,j\oplus j'}$, the pair $(i_0,j_0)$ contributes $\lambda_{i_0,j_0}\mu_{c_0',c_1}$, as $i_0+i'=c_0'<p$ and $j_0+j'=c_1<p$. Every other pair $(i,j)$ with $\lambda_{i,j}\ne 0$ has smaller index; if $j<j_0$, then $j\oplus j'=j+j'<c_1$, and if $j=j_0$ and $i<i_0$, then $i\oplus i'=i+i'<c_0'$. Hence the coefficient of $\mu_{c_0',c_1}$ in $s\ast\mu_{i',j'}$ is $\lambda_{i_0,j_0}\ne 0$. As $c_0'+pc_1\ge c$, this shows that $s\ast\mu_{i',j'}\notin\mathcal B_c$, a contradiction.
\end{proof}

\begin{proposition}\label{prop:Mforced}
Let $\bar T\le P_3$ with $\bar TH_2=P_3$ and $\bar T\cap H_2=J^c(\Omega_2)$, where $c<p^2$. Then $J^{p^2-c}(\Omega_3)$ is the only $\bar T$-submodule of $H_3=\FF_p\Omega_3$ of dimension $p^3-p^2+c$. In particular, if $T\le P_4$ has order $p^{x_4}$ and image $\bar T$ in $P_4/H_3=P_3$, then $T\cap H_3=J^{p^2-c}(\Omega_3)$.
\end{proposition}

\begin{proof}
Let $\rho\colon\FF_p\Omega_3\to\FF_p\Omega_2$ be the block sum map, let $\rho^*\colon\FF_p\Omega_2\to\FF_p\Omega_3$ be given by $e_i\mapsto\sum_{\omega\in\beta_i}e_\omega$, so that $\langle\rho^*(a),u\rangle=\langle a,\rho(u)\rangle$, and let $\mathcal C=\rho^*(\FF_p\Omega_2)$ be the space of vectors that are constant on every minimal block. The group $P_3$ acts on $\mathcal C$ through $P_3/H_2\cong P_2$, and $\rho^*$ is an isomorphism of $P_3$-modules onto $\mathcal C$; as $\bar T$ maps onto $P_3/H_2$, the $\bar T$-submodules of $\mathcal C$ are the spaces $\rho^*(J^e(\Omega_2))$.

Let $M\le H_3$ be a $\bar T$-submodule of dimension $p^3-p^2+c$, and let $L=M^\perp$, a $\bar T$-submodule of dimension $p^2-c$. Suppose that $L\not\subseteq\mathcal C$. As $\bar T$ is a $p$-group, there is $v\in L\setminus\mathcal C$ such that $vt-v\in L\cap\mathcal C$ for all $t\in\bar T$. For $i\in\Omega_2$ fix $\omega_i\in\beta_i$. An element $n=\sum_in_ie_i\in J^c(\Omega_2)=\bar T\cap H_2$, that is, $n=\prod_iy_i^{n_i}$, satisfies $(vn)_{\omega_iy_i^x}=v_{\omega_iy_i^{x-n_i}}$ for all $i$ and $x$. Since $c<p^2$, the vector $\sum_ie_i$ lies in $J^c(\Omega_2)$ by Lemma~\ref{lem:uniserial}; applying the above to it, we see that $x\mapsto v_{\omega_iy_i^{x-1}}-v_{\omega_iy_i^x}$ is constant for each $i$, so there are $s_i,d_i\in\FF_p$ with $v_{\omega_iy_i^x}=s_ix+d_i$ for all $x\in\ZZ_p$. Put $s=\sum_is_ie_i\in\FF_p\Omega_2$. Then $vn-v=-\rho^*(n\ast s)$ for all $n\in J^c(\Omega_2)$. Now $L\cap\mathcal C=\rho^*(J^e(\Omega_2))$ for some $e$, and $e\ge c+1$ as $\dim(L\cap\mathcal C)<\dim L$. Hence $J^c(\Omega_2)\ast s\subseteq J^{c+1}(\Omega_2)$, so that $s=0$ by Lemma~\ref{lem:products}. But then $v\in\mathcal C$, a contradiction.

Thus $L\subseteq\mathcal C$, and as $\dim L=p^2-c$, we obtain $L=\rho^*(J^c(\Omega_2))$. Therefore $M=L^\perp=\rho^{-1}(J^c(\Omega_2)^\perp)=\rho^{-1}(J^{p^2-c}(\Omega_2))=J^{p^2-c}(\Omega_3)$ by Lemma~\ref{lem:blocksum}. Finally, if $T$ is as in the last assertion, then $T\cap H_3$ is normalized by $T$ and centralized by $H_3$, so it is a $\bar T$-submodule of $H_3$. Since $x_4=p^3+p+1$ and $|\bar T|=|P_3|/|H_2:J^c(\Omega_2)|=p^{p^2+p+1-c}$, its dimension is $p^3-p^2+c$.
\end{proof}

\subsection{Regular orbits and a linear map}\label{sec:criterion}

Throughout this subsection, $T\le P_4$ is a subgroup of order $p^{x_4}$ having a largest orbit of $P_4$ on $V_4$ as a regular orbit, and $\bar T$ is its image in $P_4/H_3=P_3$. By Lemma~\ref{lem:image} and Section~\ref{sec:coords}, we have $\bar TH_2=P_3$ and $\bar T\cap H_2=J^c(\Omega_2)$ for some $c\le p^2-p-1$. We put
\[
N=J^c(\Omega_2)\le H_2,\qquad M=T\cap H_3,\qquad U=T\cap H_2H_3\qquad\text{and}\qquad X=H_3/M.
\]
By Proposition~\ref{prop:Mforced}, we have $M=J^{p^2-c}(\Omega_3)$. By Lemma~\ref{lem:blocksum}, $H_2$ acts trivially on $X$, and the block sum map $\rho$ induces an isomorphism from $X$ onto $\FF_p\Omega_2/J^{p^2-c}(\Omega_2)$, by which we identify these two spaces.

\begin{lemma}\label{lem:graph}
\begin{enumerate}[\rm(a)]
\item The group $H_2H_3/M$ is elementary abelian, and it is the direct product of $H_2M/M\cong H_2=\FF_p\Omega_2$ and $X=\FF_p\Omega_2/J^{p^2-c}(\Omega_2)$. Conjugation induces an action of $P_4/H_2H_3\cong P_2$ on $H_2H_3/M$, and $P_2$ acts on both factors by permuting $\Omega_2$.
\item There is a unique map $\varphi\colon N\to X$ such that, in additive notation, $U/M=\{n+\varphi(n):\ n\in N\}$; here we regard $n\in N\le H_2$ as an element of $H_2M/M$. The map $\varphi$ is linear, and $\varphi(ng)=\varphi(n)g$ for all $n\in N$ and $g\in P_2$.
\end{enumerate}
\end{lemma}

\begin{proof}
(a) As $H_2$ acts trivially on $H_3/M$, we have $[H_2,H_3]\le M$, so $H_2H_3/M$ is abelian. It is generated by the elementary abelian subgroups $H_2M/M$ and $H_3/M$, which intersect trivially because $H_2\cap H_3=1$. As $H_2H_3/M$ is abelian, $H_2H_3$ acts trivially on it by conjugation, so $P_4/H_2H_3\cong P_2$ acts on it. The last assertion follows from the identifications in Section~\ref{sec:orbits} and from the fact that $\rho$ commutes with the action of $P_3$.

(b) Since $U\cap H_3=M$, the subgroup $U/M$ meets $X$ trivially. An element of $T$ lies in $H_2H_3$ if and only if its image in $P_3$ lies in $H_2$, so the image of $U$ in $H_2H_3/H_3\cong H_2$ is $\bar T\cap H_2=N$. Hence for every $n\in N$ there is exactly one $x\in X$ with $n+x\in U/M$, and we put $\varphi(n)=x$. As $U/M$ is a subgroup, $\varphi$ is additive, hence linear. Finally, let $g\in P_2$. As $\bar TH_2=P_3$, we have $TH_2H_3=P_4$, so there is $t\in T$ with $tH_2H_3=gH_2H_3$. Conjugation by $t$ maps $U$ and $M$ to themselves, so it maps $U/M$ to itself, and by (a) it maps $n+\varphi(n)$ to $ng+\varphi(n)g$. Hence $\varphi(ng)=\varphi(n)g$.
\end{proof}

\begin{proposition}\label{prop:criterion}
Let $\Delta'\subseteq\Omega_2$ lie in a regular orbit of $P_2$ on $\Pow(\Omega_2)$, put $W_{\Delta'}=N\cap\FF_p\Delta'$, and let $X_{\Delta'}$ be the image of $\FF_p(\Omega_2\setminus\Delta')$ in $X=\FF_p\Omega_2/J^{p^2-c}(\Omega_2)$. Then $O_{\Delta'}$ is a regular orbit of $T$ if and only if
\begin{enumerate}[\rm(i)]
\item $N+\FF_p\Delta'=\FF_p\Omega_2$, and
\item $\varphi(n)\notin X_{\Delta'}$ for every $0\ne n\in W_{\Delta'}$.
\end{enumerate}
\end{proposition}

\begin{proof}
By Theorem~\ref{thm:largest}(d), we have $C:=\C_{P_4}(v_{\Delta'})=Y_{\Delta'}Z_{\Delta'}\le H_2H_3$, where $Y_{\Delta'}=\FF_p\Delta'\le H_2$ and $Z_{\Delta'}=\FF_p\Sigma_{\Delta'}\le H_3$ with $\Sigma_{\Delta'}=\{\sigma_i:\ i\notin\Delta'\}$. In particular $T\cap C=U\cap C$, so by Lemma~\ref{lem:regular}, $O_{\Delta'}$ is a regular orbit of $T$ if and only if $U\cap C=1$. Every element of $U$ can be written uniquely as $nh$ with $n\in N$ and $h\in H_3$, and then $hM=\varphi(n)$ by Lemma~\ref{lem:graph}(b); this element lies in $C$ if and only if $n\in\FF_p\Delta'$ and $h\in\FF_p\Sigma_{\Delta'}$. As $\rho$ maps $\FF_p\Sigma_{\Delta'}$ isomorphically onto $\FF_p(\Omega_2\setminus\Delta')$, the image of $\FF_p\Sigma_{\Delta'}$ in $X$ is $X_{\Delta'}$. Hence $U\cap C=1$ if and only if $M\cap\FF_p\Sigma_{\Delta'}=1$ and (ii) holds. Moreover, as $M=\rho^{-1}(J^{p^2-c}(\Omega_2))$, we have $M\cap\FF_p\Sigma_{\Delta'}=1$ if and only if $J^{p^2-c}(\Omega_2)\cap\FF_p(\Omega_2\setminus\Delta')=0$, and taking orthogonal complements (as in the proof of Lemma~\ref{lem:goodcomplement}), this holds if and only if (i) holds.
\end{proof}

To make condition (ii) computable, we use the standard bilinear form.

\begin{proposition}\label{prop:matrix}
For $n,n'\in N$ put $\Phi(n,n')=\langle a,n'\rangle$, where $a\in\FF_p\Omega_2$ is any vector with $a+J^{p^2-c}(\Omega_2)=\varphi(n)$. Then $\Phi\colon N\times N\to\FF_p$ is well defined and bilinear, and $\Phi(ng,n'g)=\Phi(n,n')$ for all $n,n'\in N$ and $g\in P_2$. Moreover, let $\Delta'$ be as in Proposition~\ref{prop:criterion} such that \textup{(i)} holds, and let $w_1,\dots,w_r$ be a basis of $W_{\Delta'}$. Then $O_{\Delta'}$ is a regular orbit of $T$ if and only if the $r\times r$ matrix $(\Phi(w_i,w_j))_{i,j}$ is invertible.
\end{proposition}

\begin{proof}
As $J^{p^2-c}(\Omega_2)=N^\perp$, the map $\Phi$ is well defined, and it is clearly bilinear; the invariance follows from Lemma~\ref{lem:graph}(b), as $\langle ag,n'g\rangle=\langle a,n'\rangle$. Assigning to $a+J^{p^2-c}(\Omega_2)\in X$ the linear map $N\to\FF_p$, $n'\mapsto\langle a,n'\rangle$, we obtain an injective linear map from $X$ to the dual space of $N$, which is bijective as $\dim X=p^2-c=\dim N$. Every element of $X_{\Delta'}$ is mapped to a linear map vanishing on $W_{\Delta'}$, as the vectors in $\FF_p(\Omega_2\setminus\Delta')$ are orthogonal to $\FF_p\Delta'$. By (i), we have $\dim W_{\Delta'}=|\Delta'|-c$, and by the proof of Proposition~\ref{prop:criterion}, $J^{p^2-c}(\Omega_2)\cap\FF_p(\Omega_2\setminus\Delta')=0$, so that $\dim X_{\Delta'}=p^2-|\Delta'|=\dim N-\dim W_{\Delta'}$. This is the dimension of the space of all linear maps $N\to\FF_p$ that vanish on $W_{\Delta'}$, so $X_{\Delta'}$ corresponds exactly to this space. Therefore condition (ii) of Proposition~\ref{prop:criterion} says that for every $0\ne n\in W_{\Delta'}$ there is $n'\in W_{\Delta'}$ with $\Phi(n,n')\ne 0$. Writing $n=\sum_i\lambda_iw_i$, this means that no nontrivial linear combination of the rows of the matrix $(\Phi(w_i,w_j))_{i,j}$ vanishes, that is, that this matrix is invertible.
\end{proof}

\subsection{The shape of the matrices}\label{sec:forms}

Let $G$ be a group of permutations of a finite set $\Omega$, and let $J\subseteq\FF_p\Omega$ be a subspace with $Jg=J$ for all $g\in G$. We call a bilinear map $\Phi\colon J\times J\to\FF_p$ \emph{$G$-invariant} if $\Phi(ug,wg)=\Phi(u,w)$ for all $u,w\in J$ and $g\in G$; for example, the map $\Phi$ of Proposition~\ref{prop:matrix} is $P_2$-invariant. If $\Phi(u,w)=u\Psi w^{\mathsf T}$ for a matrix $\Psi$ indexed by $\Omega$, where $u$ and $w$ are regarded as row vectors, then $\Phi$ is $G$-invariant if $\Psi$ commutes with the permutation matrices of the elements of $G$. For $G=P_1$ acting on $\Omega_1=\ZZ_p$, the matrices commuting with the permutation matrix of $x\mapsto x+1$ are the circulant matrices $(a(x'-x))_{x,x'\in\ZZ_p}$ with $a\colon\ZZ_p\to\FF_p$. The following two lemmas show that, conversely, the $P_2$-invariant bilinear maps on $J^c(\Omega_2)$ with $c\le p-1$ are given by matrices built from circulant matrices.

\begin{lemma}\label{lem:circrestr}
Let $\kappa\colon J^1(\Omega_1)\times J^1(\Omega_1)\to\FF_p$ be a $P_1$-invariant bilinear map such that $\kappa(\one,e_1-e_0)=0$, where $\one=\sum_{x\in\ZZ_p}e_x$. Then there is a map $a\colon\ZZ_p\to\FF_p$ such that $\kappa(z,z')=\sum_{x,x'\in\ZZ_p}a(x'-x)z_xz'_{x'}$ for all $z,z'\in J^1(\Omega_1)$.
\end{lemma}

\begin{proof}
For $x\in\ZZ_p$ put $\delta_x=e_{x+1}-e_x$; these vectors span $J^1(\Omega_1)$, and $P_1$ permutes them regularly. Hence $\kappa(\delta_x,\delta_{x'})=k(x'-x)$, where $k(d)=\kappa(\delta_0,\delta_d)$. As $\sum_d\delta_d=0$, we have $\sum_dk(d)=0$. Moreover, $\sum_x(x+1)\delta_x=-\one$, so that $0=\kappa(\one,\delta_0)=-\sum_x(x+1)k(-x)=\sum_dd\,k(d)-\sum_dk(d)=\sum_dd\,k(d)$. Hence the polynomial $f=\sum_{d=0}^{p-1}k(d)t^d$ satisfies $f(1)=f'(1)=0$, so $f$ lies in the ideal $(t-1)^2R$ of $R=\FF_p[t]/(t-1)^p=\FF_p[t]/(t^p-1)$. As $t$ is a unit of $R$, there is $A=\sum_{d\in\ZZ_p}a(d)t^d\in R$ with $f=-(t-2+t^{-1})A$, that is, $k(d)=2a(d)-a(d-1)-a(d+1)$ for all $d\in\ZZ_p$. The $P_1$-invariant bilinear map $\kappa_a(z,z')=\sum_{x,x'}a(x'-x)z_xz'_{x'}$ satisfies $\kappa_a(\delta_0,\delta_d)=2a(d)-a(d+1)-a(d-1)=k(d)$. Hence $\kappa$ and $\kappa_a$ agree on all pairs $(\delta_x,\delta_{x'})$, and therefore on $J^1(\Omega_1)$.
\end{proof}

\begin{lemma}\label{lem:forms}
Let $0\le c\le p$, and let $\Phi\colon J\times J\to\FF_p$ be a $P_2$-invariant bilinear map, where $J=J^c(\Omega_2)$.
\begin{enumerate}[\rm(a)]
\item If $c\le p-1$, then there are maps $\psi\colon\ZZ_p\to\FF_p$ and $\eta\colon\ZZ_p\times\ZZ_p\to\FF_p$ such that for all $u,w\in J$ we have
\[
\Phi(u,w)=\sum_{b\in\ZZ_p}\ \sum_{x,x'\in\ZZ_p}\psi(x'-x)\,u_{(b,x)}w_{(b,x')}+\sum_{b,b'\in\ZZ_p}\eta(b,b')\,\varepsilon_b(u)\varepsilon_{b'}(w).
\]
\item If $c=p$ and $b\ne b'$, then there is $\gamma\in\FF_p$ such that $\Phi(u,w)=\gamma\,\ell_b(u)\ell_{b'}(w)$ for all $u\in K_b$ and $w\in K_{b'}$.
\end{enumerate}
\end{lemma}

\begin{proof}
(b) Let $u\in K_b$ and $w\in K_{b'}$. As $uh_{b'}=u$, invariance gives $\Phi(u,wh_{b'}-w)=0$, so by \eqref{eq:Kb} the linear form $\Phi(u,\cdot)$ on $K_{b'}$ vanishes on the kernel of $\ell_{b'}$ and is therefore a multiple of $\ell_{b'}$. Similarly, $\Phi(\cdot,w)$ is a multiple of $\ell_b$ on $K_b$. Choosing $u_1\in K_b$ and $w_1\in K_{b'}$ with $\ell_b(u_1)=\ell_{b'}(w_1)=1$, we obtain $\Phi(u,w)=\ell_b(u)\Phi(u_1,w)=\ell_b(u)\ell_{b'}(w)\Phi(u_1,w_1)$.

(a) Let $c\le p-1$. Then $\one=\sum_be_b$ lies in $J^c(\Omega_1)$ by Lemma~\ref{lem:uniserial}, so for every $b$ there are vectors $w\in J$ with $\varepsilon_b(w)=1$; moreover, $K\subseteq J$ by \eqref{eq:Jc}. We proceed in several steps.

\emph{Step 1. If $b\in\ZZ_p$ and $u\in J$ satisfies $u_{(b,x)}=0$ for all $x$, then $\Phi(u,z)=\Phi(z,u)=0$ for all $z\in K_b$.} Indeed, $uh_b=u$, so $\Phi(u,wh_b-w)=0$ for all $w\in J$. The vectors $wh_b-w$ lie in $K_b$ and satisfy $\ell_b(wh_b-w)=\varepsilon_b(w)$; by \eqref{eq:Kb} and the remark above, they span $K_b$. The second assertion follows in the same way.

\emph{Step 2.} For $u\in\FF_p\Omega_2$ and $b\in\ZZ_p$ let $k_b(u)=\sum_xu_{(b,x)}e_{(b,x)}-\varepsilon_b(u)e_{(b,0)}\in K_b$, and for $v\in\FF_p\Omega_1$ let $s(v)=\sum_bv_be_{(b,0)}$. Then $u=\sum_bk_b(u)+s(\rho_2(u))$, and $s(\rho_2(u))\in J$ if $u\in J$. By Step 1, $\Phi(k_b(u),k_{b'}(w))=0$ for $b\ne b'$. For $z\in K_b$, the linear form $v\mapsto\Phi(z,s(v))$ on $J^c(\Omega_1)$ vanishes on $\{v\in J^c(\Omega_1):\ v_b=0\}$ by Step 1, and this is a hyperplane of $J^c(\Omega_1)$ as $\one\in J^c(\Omega_1)$; hence $\Phi(z,s(v))=\alpha'_b(z)v_b$ for a linear form $\alpha'_b$ on $K_b$. Similarly, $\Phi(s(v),z)=v_b\alpha_b(z)$. Let $\kappa_b$ be the restriction of $\Phi$ to $K_b\times K_b$, and put $\Phi_0(v,v')=\Phi(s(v),s(v'))$ for $v,v'\in J^c(\Omega_1)$. Then, for all $u,w\in J$, we have
\begin{multline}\label{eq:Phidecomp}
\Phi(u,w)=\sum_b\Bigl(\kappa_b\bigl(k_b(u),k_b(w)\bigr)+\alpha'_b\bigl(k_b(u)\bigr)\varepsilon_b(w)+\varepsilon_b(u)\alpha_b\bigl(k_b(w)\bigr)\Bigr)\\
+\Phi_0\bigl(\rho_2(u),\rho_2(w)\bigr).
\end{multline}

\emph{Step 3.} Fix $b$, put $\delta_b=e_{(b,1)}-e_{(b,0)}$, choose $w\in J$ with $\varepsilon_b(w)=1$, and let $s_w=s(\rho_2(w))$. Then $s_wh_b=s_w+\delta_b$, so for $z\in K_b$ invariance under $h_b$ yields
\begin{equation}\label{eq:alphashift}
\alpha'_b(z)=\Phi(z,s_w)=\Phi(zh_b,s_w+\delta_b)=\alpha'_b(zh_b)+\kappa_b(zh_b,\delta_b).
\end{equation}
For $z=\one_b=\sum_xe_{(b,x)}\in K_b$ we have $zh_b=z$, hence $\kappa_b(\one_b,\delta_b)=0$. By Lemma~\ref{lem:circrestr}, applied via the bijection $B_b\to\ZZ_p$, $(b,x)\mapsto x$, there is $a_b\colon\ZZ_p\to\FF_p$ such that $\kappa_b$ is the restriction to $K_b\times K_b$ of the bilinear map $C_b(u,w)=\sum_{x,x'}a_b(x'-x)u_{(b,x)}w_{(b,x')}$ on $\FF_p\Omega_2$. As $\hat h$ maps $(b,x)$ to $(b+1,x)$, invariance under $\hat h$ shows that we may choose $a_b=a$ independently of $b$, and that $\alpha_{b+1}(z\hat h)=\alpha_b(z)$ and $\alpha'_{b+1}(z\hat h)=\alpha'_b(z)$ for $z\in K_b$.

\emph{Step 4.} The bilinear map $C_b$ is invariant under $h_b$, so $z\mapsto C_b(z,e_{(b,0)})$ also satisfies \eqref{eq:alphashift}: $C_b(z,e_{(b,0)})=C_b(zh_b,e_{(b,1)})=C_b(zh_b,e_{(b,0)})+\kappa_b(zh_b,\delta_b)$. Hence $\theta_b(z)=\alpha'_b(z)-C_b(z,e_{(b,0)})$ satisfies $\theta_b(z)=\theta_b(zh_b)$ for $z\in K_b$, so $\theta_b$ vanishes on $\{zh_b-z:\ z\in K_b\}$, and $\theta_b=\lambda'_b\ell_b$ on $K_b$ for some $\lambda'_b\in\FF_p$ by \eqref{eq:Kb}. Similarly, $\alpha_b(z)=C_b(e_{(b,0)},z)+\lambda_b\ell_b(z)$ for some $\lambda_b\in\FF_p$. By Step 3, $\lambda_b=\lambda$ and $\lambda'_b=\lambda'$ do not depend on $b$.

\emph{Step 5.} Choose $u\in J$ with $\varepsilon_b(u)=1$, let $s_u=s(\rho_2(u))$, and let $w$ and $s_w$ be as in Step 3. Invariance gives $\Phi(s_u+\delta_b,s_w+\delta_b)=\Phi(s_u,s_w)$, that is, $\alpha_b(\delta_b)+\alpha'_b(\delta_b)+\kappa_b(\delta_b,\delta_b)=0$. The corresponding identity $C_b(e_{(b,0)},\delta_b)+C_b(\delta_b,e_{(b,0)})+C_b(\delta_b,\delta_b)=0$ follows from $C_b(e_{(b,1)},e_{(b,1)})=C_b(e_{(b,0)},e_{(b,0)})$. Subtracting the two identities and using $\ell_b(\delta_b)=1$, we obtain $\lambda+\lambda'=0$.

\emph{Step 6.} For $u,w\in\FF_p\Omega_2$ we have $\ell_b(k_b(u))=\ell_b(u)$ and
\begin{multline*}
C_b(u,w)=C_b\bigl(k_b(u),k_b(w)\bigr)+C_b\bigl(k_b(u),e_{(b,0)}\bigr)\varepsilon_b(w)\\
+\varepsilon_b(u)C_b\bigl(e_{(b,0)},k_b(w)\bigr)+a(0)\varepsilon_b(u)\varepsilon_b(w).
\end{multline*}
Substituting Steps 3--5 into \eqref{eq:Phidecomp}, we obtain for $u,w\in J$
\begin{multline*}
\Phi(u,w)=\sum_b\Bigl(C_b(u,w)+\lambda\bigl(\varepsilon_b(u)\ell_b(w)-\ell_b(u)\varepsilon_b(w)\bigr)-a(0)\varepsilon_b(u)\varepsilon_b(w)\Bigr)\\
+\Phi_0\bigl(\rho_2(u),\rho_2(w)\bigr).
\end{multline*}
Since $\varepsilon_b(u)\ell_b(w)-\ell_b(u)\varepsilon_b(w)=\sum_{x,x'}(x'-x)u_{(b,x)}w_{(b,x')}$, the assertion follows with $\psi(d)=a(d)+\lambda d$ and any $\eta$ such that $\sum_{b,b'}\eta(b,b')v_bv'_{b'}=\Phi_0(v,v')-a(0)\sum_bv_bv'_b$ for all $v,v'\in J^c(\Omega_1)$.
\end{proof}

\subsection{One block at a time}\label{sec:oneblock}

In this subsection, $\FF$ is an arbitrary field. The following lemma describes how the determinant of a matrix depends on a block of rows and columns that is coupled to the remaining rows and columns only through a fixed vector.

\begin{lemma}\label{lem:oneblock}
Let $G$ be a square matrix over $\FF$ whose rows and columns are indexed by a finite set $I_0\sqcup I_1$ with $I_0\ne\emptyset$. Let $u\in\FF^{I_0}$ and $i_0\in I_0$ with $u_{i_0}\ne 0$, and assume that there are $y,z\in\FF^{I_1}$ such that $G_{ij}=u_iy_j$ and $G_{ji}=z_ju_i$ for all $i\in I_0$ and $j\in I_1$. Let $A$ and $R$ be the submatrices of $G$ with rows and columns indexed by $I_0$ and by $I_1$, respectively, and regard $u,y,z$ as column vectors. Then:
\begin{enumerate}[\rm(a)]
\item There are $\alpha,\beta\in\FF$ such that $\det(A+t\,uu^{\mathsf T})=\alpha+\beta t$ for all $t\in\FF$.
\item $\det G=d\alpha+e\beta$, where $d=\det R$ and $e=\det\begin{pmatrix}0&y^{\mathsf T}\\ z&R\end{pmatrix}$; here $d=1$ and $e=0$ if $I_1=\emptyset$.
\end{enumerate}
\end{lemma}

\begin{proof}
For each $i\in I_0\setminus\{i_0\}$, subtract $u_i/u_{i_0}$ times the row $i_0$ from the row $i$, and then $u_i/u_{i_0}$ times the column $i_0$ from the column $i$. These operations do not change determinants. They transform $A$ into a matrix $\tilde A$, the matrix $A+t\,uu^{\mathsf T}$ into $\tilde A+tu_{i_0}^2E$, where $E$ is the matrix unit at the position $(i_0,i_0)$, and $G$ into a matrix $\tilde G$ in which the rows and the columns indexed by $I_0\setminus\{i_0\}$ have zero entries outside $I_0\times I_0$, while the row and the column $i_0$ have the entries $u_{i_0}y_j$ and $u_{i_0}z_j$, respectively, at $j\in I_1$. Let $B$ be the submatrix of $\tilde A$ with rows and columns indexed by $I_0\setminus\{i_0\}$, where $\det B=1$ if $I_0=\{i_0\}$. Expanding along the row $i_0$, we obtain $\det(\tilde A+tu_{i_0}^2E)=\det\tilde A+tu_{i_0}^2\det B$, which proves (a) with $\alpha=\det A$ and $\beta=u_{i_0}^2\det B$.

For (b), we write the row $i_0$ of $\tilde G$ as the sum of its part inside $I_0$ and its part inside $I_1$; then $\det\tilde G=\det G_1+\det G_2$, where the row $i_0$ of $G_1$ (resp.\ $G_2$) vanishes outside $I_0$ (resp.\ inside $I_0$). The matrix $G_1$ is block triangular with diagonal blocks $\tilde A$ and $R$, so $\det G_1=\alpha d$. Splitting the column $i_0$ of $G_2$ in the same way, we obtain $\det G_2=\det G_3+\det G_4$, where the column $i_0$ of $G_3$ (resp.\ $G_4$) vanishes outside $I_0$ (resp.\ inside $I_0$). In $G_3$, the $|I_1|+1$ rows indexed by $\{i_0\}\cup I_1$ have zero entries outside the $|I_1|$ columns indexed by $I_1$, so $\det G_3=0$. In $G_4$, the rows and the columns indexed by $I_0\setminus\{i_0\}$ have zero entries outside $B$, so $\det G_4=\det B\cdot\det\begin{pmatrix}0&u_{i_0}y^{\mathsf T}\\ u_{i_0}z&R\end{pmatrix}=\beta e$.
\end{proof}

In the situation of Lemma~\ref{lem:oneblock}(a), we call the pair $(A,u)$ \emph{dead} if $\alpha=\beta=0$, that is, if $A+t\,uu^{\mathsf T}$ is singular for every $t\in\FF$, and \emph{live} otherwise. If $(A,u)$ is live, its \emph{forbidden value} $f(A,u)\in\FF\cup\{\infty\}$ is the root $-\alpha/\beta$ of $\alpha+\beta t$ if $\beta\ne 0$, and $f(A,u)=\infty$ if $\beta=0$. Informally, the following corollary says that, from the point of view of the block $I_0$, the remaining rows and columns of $G$ act like adding $\theta\,uu^{\mathsf T}$ to $A$ for a suitable $\theta$, so that $\det G\ne 0$ unless $\theta$ is the forbidden value.

\begin{corollary}\label{cor:forbidden}
In the situation of Lemma~\ref{lem:oneblock}, put $\theta=e/d$ if $d\ne 0$, $\theta=\infty$ if $d=0\ne e$, and $\theta=\ast$ if $d=e=0$. Then $\det G\ne 0$ if and only if $\theta\ne\ast$, the pair $(A,u)$ is live, and $f(A,u)\ne\theta$.
\end{corollary}

\begin{proof}
If $d\ne 0$, then $\det G=d(\alpha+\theta\beta)$, which vanishes if and only if $(A,u)$ is dead or $f(A,u)=\theta$. If $d=0\ne e$, then $\det G=e\beta$, which vanishes if and only if $\beta=0$, that is, if and only if $(A,u)$ is dead or $f(A,u)=\infty=\theta$. If $d=e=0$, then $\det G=0$.
\end{proof}

\subsection{Small circulant matrices}\label{sec:circulants}

For $\psi\colon\ZZ_p\to\FF_p$ and $\Gamma\subseteq\ZZ_p$ let $C_\Gamma(\psi)=(\psi(x'-x))_{x,x'\in\Gamma}$, with rows and columns indexed by $\Gamma$, and let $\one$ be the all-one column vector indexed by $\Gamma$. If $\Gamma'=\Gamma+1$, then $C_{\Gamma'}(\psi)$ arises from $C_\Gamma(\psi)$ by renaming the indices. So whether the pair $(C_\Gamma(\psi),\one)$ is live, and its forbidden value, only depend on the type $\tau$ of $\Gamma$; we call $\tau$ \emph{$\psi$-live} or \emph{$\psi$-dead} accordingly, and for a $\psi$-live type $\tau$ we denote this forbidden value by $f_\psi(\tau)$. Note that $C_\Gamma(\psi)+t\,\one\one^{\mathsf T}=C_\Gamma(\psi+t)$, where $\psi+t$ is the map $x\mapsto\psi(x)+t$.

\begin{lemma}\label{lem:circulant}
Let $p\ge 5$ and $\psi\colon\ZZ_p\to\FF_p$. If $\psi$ is constant, then the type of the $1$-element sets is the only $\psi$-live type. Otherwise, there are at least three $\psi$-live types.
\end{lemma}

\begin{proof}
Let $t\in\FF_p$ and $s=\psi(0)+t$. For $\Gamma=\{0\}$ we have $\det C_\Gamma(\psi+t)=s$, so the type of the $1$-element sets is live. For $0\ne d\in\ZZ_p$ we have
\begin{align*}
\det C_{\{0,d\}}(\psi+t)&=s^2-\bigl(\psi(d)+t\bigr)\bigl(\psi(-d)+t\bigr)\\
&=\bigl(\psi(0)^2-\psi(d)\psi(-d)\bigr)+\bigl(2\psi(0)-\psi(d)-\psi(-d)\bigr)t,
\end{align*}
which vanishes for all $t$ if and only if $\psi(d)$ and $\psi(-d)$ are the roots of $(X-\psi(0))^2$, that is, if and only if $\psi(d)=\psi(-d)=\psi(0)$. The types of $2$-element sets are those of the sets $\{0,d\}$ with $1\le d\le(p-1)/2$, so there are $(p-1)/2\ge 2$ of them.

If $\psi$ is constant, then $C_\Gamma(\psi+t)=s\,\one\one^{\mathsf T}$ has rank at most one, so it is singular for every $t$ whenever $|\Gamma|\ge 2$. Hence only the $1$-element sets are live.

Now let $\psi$ be non-constant. Then some type of $2$-element sets is live. If two of them are live, we are done; so we may assume that only the type of $\{0,d\}$ is live, and then $\psi(e)=\psi(0)$ for all $e\notin\{0,d,-d\}$. The map $x\mapsto d^{-1}x$ permutes the types, and $C_{d^{-1}\Gamma}(\psi_d)$ arises from $C_\Gamma(\psi)$ by renaming the indices, where $\psi_d(x)=\psi(dx)$. So, replacing $\psi$ by $\psi_d$, we may assume that $d=1$. Write $\psi(1)=\psi(0)+u$ and $\psi(-1)=\psi(0)+u'$, so that $(u,u')\ne(0,0)$. As $p\ge 5$, two of the elements $0,1,2,3$ of $\ZZ_p$ differ by $\pm1$ only if they are consecutive integers. Hence, with the rows and columns ordered as written,
\[
C_{\{0,1,3\}}(\psi+t)=\begin{pmatrix}s&s+u&s\\ s+u'&s&s\\ s&s&s\end{pmatrix},\qquad
C_{\{0,1,2\}}(\psi+t)=\begin{pmatrix}s&s+u&s\\ s+u'&s&s+u\\ s&s+u'&s\end{pmatrix}.
\]
Subtracting the last row from the other two rows, we obtain $\det C_{\{0,1,3\}}(\psi+t)=-uu's$, and $\det C_{\{0,1,2\}}(\psi+t)=u^2s$ if $u'=0$, and $=u'^2s$ if $u=0$. Hence the type of $\{0,1,3\}$ is live if $uu'\ne 0$, and the type of $\{0,1,2\}$ is live if $uu'=0$. In each case we have found a live type of $3$-element sets, in addition to the two live types found above.
\end{proof}

\subsection{Proof of Theorem~\ref{thmE}}\label{sec:proofE}

The following combinatorial lemma is the heart of the proof.

\begin{lemma}\label{lem:blockrepl}
Suppose that every type is either dead or live, and that every live type $\tau$ has a forbidden value $f(\tau)\in\FF_p\cup\{\infty\}$. Let $\mathcal R$ be a set of non-constant $p$-tuples of types that is closed under cyclic shifts. Assume that for every non-constant tuple $\mathbf t$ and every $b\in\ZZ_p$ there is $\theta_b(\mathbf t)\in\FF_p\cup\{\infty,\ast\}$, depending only on $b$ and on the entries $t_{b'}$ with $b'\ne b$, such that $\mathbf t\in\mathcal R$ if and only if $\theta_b(\mathbf t)\ne\ast$, the type $t_b$ is live, and $f(t_b)\ne\theta_b(\mathbf t)$. If the number of live types is different from $2$, then $\mathcal R$ is not a single cyclic-shift class.
\end{lemma}

\begin{proof}
Assume that $\mathcal R$ is the cyclic-shift class of $\mathbf t^*=(t^*_b)$. Every non-constant tuple has two different entries, so if at most one type is live, then $\mathcal R=\emptyset$, a contradiction. Hence there are at least three live types. For $b\in\ZZ_p$ put $\theta_b=\theta_b(\mathbf t^*)$; then $\theta_b\ne\ast$, the type $t^*_b$ is live, and $f(t^*_b)\ne\theta_b$. Call a type $\tau$ \emph{exceptional} for $b$ if $t^*_{b'}=\tau$ for all $b'\ne b$. Let $\tau\ne t^*_b$ be a live type that is not exceptional for $b$, and let $\mathbf t'$ arise from $\mathbf t^*$ by replacing $t^*_b$ by $\tau$. Then $\mathbf t'$ is non-constant and not a cyclic shift of $\mathbf t^*$, as the multisets of entries differ, so $\mathbf t'\notin\mathcal R$; as $\theta_b(\mathbf t')=\theta_b$, this means that $f(\tau)=\theta_b$. Thus:
\begin{equation}\label{eq:star}
f(\tau)=\theta_b\quad\text{for every live type }\tau\ne t^*_b\text{ that is not exceptional for }b.
\end{equation}
If some $b$ has an exceptional type, then all entries of $\mathbf t^*$ except $t^*_b$ coincide.

Assume first that no $b$ has an exceptional type. Choose $b\ne b'$ with $t^*_b\ne t^*_{b'}$ and a live type $\tau\notin\{t^*_b,t^*_{b'}\}$. By \eqref{eq:star}, $\theta_b=f(\tau)=\theta_{b'}$, and also $f(t^*_{b'})=\theta_b=\theta_{b'}$, a contradiction.

Otherwise there are $b_0\in\ZZ_p$ and types $\tau_1\ne\tau_2$ with $t^*_{b_0}=\tau_1$ and $t^*_b=\tau_2$ for all $b\ne b_0$. Then $\tau_2$ is the only exceptional type for $b_0$, and no type is exceptional for $b\ne b_0$, since the entries $t^*_{b'}$ with $b'\ne b$ include $\tau_1$ and, as $p\ge 3$, also $\tau_2$. Fix $b\ne b_0$ and a live type $\tau_3\notin\{\tau_1,\tau_2\}$. By \eqref{eq:star}, applied to $b$ and to $b_0$, we have $f(\tau_1)=\theta_b=f(\tau_3)=\theta_{b_0}$, that is, $f(t^*_{b_0})=\theta_{b_0}$, a contradiction.
\end{proof}

\begin{lemma}\label{lem:twogood}
Let $p\ge 5$ and $p+1\le c\le p^2-p-1$. Then there are non-constant tuples $\mathbf t$ and $\mathbf t'$ of types, which are not cyclic shifts of each other, such that $\Delta(\mathbf t)$ and $\Delta(\mathbf t')$ are balanced subsets of $\Omega_2$ of size $c$.
\end{lemma}

\begin{proof}
Write $c=pq+r$ with $0\le r<p$. Then $1\le q\le p-2$, and $q\ge 2$ if $r=0$. Let $s=q$ if $q\ge 2$, and $s=q+1=2$ otherwise (in which case $r\ge 1$). Then $2\le s\le p-2$, so there are $\binom ps/p\ge(p-1)/2\ge 2$ types of $s$-element sets; let $\tau_1\ne\tau_2$ be two of them. If $r\ge 1$, let $\mathbf t$ have $r$ entries of size $q+1$ and $p-r$ entries of size $q$ (meaning types of sets of these sizes), where all entries of size $s$ equal $\tau_1$, and let $\mathbf t'$ arise from $\mathbf t$ by replacing one entry $\tau_1$ by $\tau_2$; both tuples are non-constant, as two sizes occur. If $r=0$, let $\mathbf t$ have one entry $\tau_2$ and $p-1$ entries $\tau_1$, and let $\mathbf t'$ have two entries $\tau_2$ and $p-2$ entries $\tau_1$. In both cases $\Delta(\mathbf t)$ and $\Delta(\mathbf t')$ are balanced of size $c$, and $\mathbf t'$ is not a cyclic shift of $\mathbf t$, as the multisets of entries differ.
\end{proof}

\begin{proof}[Proof of Theorem~\ref{thmE}]
Let $p\ge 7$. By Corollary~\ref{cor:descent}, it suffices to show that no subgroup $T\le P_4$ has a largest orbit $O$ of $P_4$ on $V_4$ as its unique regular orbit. Suppose that $T$ is such a subgroup. Then $|T|=|O|=p^{x_4}$, and we use the notation of Section~\ref{sec:criterion}. For a non-constant tuple $\mathbf t$ we put $W_{\mathbf t}=N\cap\FF_p\Delta(\mathbf t)$. By Section~\ref{sec:coords}, the regular orbits of $T$ are among the orbits $O_{\Delta(\mathbf t)}$. Let $\mathcal R$ be the set of non-constant tuples $\mathbf t$ for which $O_{\Delta(\mathbf t)}$ is a regular orbit of $T$; then $\mathcal R$ is closed under cyclic shifts, and by our assumption it is a single cyclic-shift class.

\emph{Case 1: $c\ge p+1$.} Let $\mathbf t$ and $\mathbf t'$ be as in Lemma~\ref{lem:twogood}. By Proposition~\ref{prop:balancedgood}, $\Delta(\mathbf t)$ is good, so $W_{\mathbf t}=0$, and condition (i) of Proposition~\ref{prop:criterion} holds as $|\Delta(\mathbf t)|=c$, while condition (ii) holds trivially. Hence $\mathbf t\in\mathcal R$, and similarly $\mathbf t'\in\mathcal R$, a contradiction.

\emph{Case 2: $c\le p$.} Every set $\Gamma_\tau$ is nonempty, so $\rho_2(\FF_p\Delta(\mathbf t))=\FF_p\Omega_1$, and condition (i) of Proposition~\ref{prop:criterion} holds for all $\mathbf t$ by \eqref{eq:Jc}. Let $\Phi$ be as in Proposition~\ref{prop:matrix}. Then $\mathbf t\in\mathcal R$ if and only if the matrix $(\Phi(w_i,w_j))_{i,j}$ is invertible for some (equivalently, every) basis $(w_i)$ of $W_{\mathbf t}$.

\emph{Case 2a: $c\le p-1$.} Let $\psi$ and $\eta$ be as in Lemma~\ref{lem:forms}(a), and let $\Psi$ be the matrix indexed by $\Omega_2$ with entries $\Psi_{(b,x),(b',x')}=\delta_{bb'}\psi(x'-x)+\eta(b,b')$, where $\delta_{bb'}$ is the Kronecker delta; then $\Phi(u,w)=u\Psi w^{\mathsf T}$ for $u,w\in N$, regarded as row vectors. Choose vectors $f_1,\dots,f_c\in\FF_p\Omega_1$ such that $J^c(\Omega_1)=\{v:\ \langle v,f_j\rangle=0\text{ for all }j\}$, for instance a basis of $J^c(\Omega_1)^\perp=J^{p-c}(\Omega_1)$. By \eqref{eq:Jc}, a vector $w\in\FF_p\Delta(\mathbf t)$ lies in $W_{\mathbf t}$ if and only if $\langle\rho_2(w),f_j\rangle=0$ for all $j$. Order $\Delta(\mathbf t)$ block by block, and let $G(\mathbf t)$ be the matrix
\[
G(\mathbf t)=\begin{pmatrix}\bigl(\Psi_{\omega\omega'}\bigr)_{\omega,\omega'\in\Delta(\mathbf t)}&F^{\mathsf T}\\ F&0\end{pmatrix},\qquad\text{where } F=\bigl((f_j)_b\bigr)_{1\le j\le c,\ (b,x)\in\Delta(\mathbf t)},
\]
so that $wF^{\mathsf T}=(\langle\rho_2(w),f_j\rangle)_j$ for $w\in\FF_p\Delta(\mathbf t)$, regarded as a row vector. Since $\rho_2(\FF_p\Delta(\mathbf t))=\FF_p\Omega_1$, the $c$ linear conditions $\langle\rho_2(w),f_j\rangle=0$ on $\FF_p\Delta(\mathbf t)$ are linearly independent. So we may choose a basis $w_1,\dots,w_r$ of $W_{\mathbf t}$ and vectors $v_1,\dots,v_c\in\FF_p\Delta(\mathbf t)$ with $\langle\rho_2(v_i),f_j\rangle=\delta_{ij}$. Let $P$ be the matrix with the rows $w_1,\dots,w_r,v_1,\dots,v_c$. Then
\[
\begin{pmatrix}P&0\\0&I\end{pmatrix}G(\mathbf t)\begin{pmatrix}P&0\\0&I\end{pmatrix}^{\mathsf T}=\begin{pmatrix}A_{11}&A_{12}&0\\ A_{21}&A_{22}&I\\ 0&I&0\end{pmatrix},
\]
where $A_{11}=(\Phi(w_i,w_j))_{i,j}$. Taking the block rows in the order $2,3,1$ and the block columns in the order $3,2,1$, we obtain the matrix
\[
\begin{pmatrix}I&A_{22}&A_{21}\\ 0&I&0\\ 0&A_{12}&A_{11}\end{pmatrix},
\]
whose determinant equals $\det A_{11}$, since the lower right part $\left(\begin{smallmatrix}I&0\\ A_{12}&A_{11}\end{smallmatrix}\right)$ is block triangular. Hence $\det G(\mathbf t)=\pm\det(P)^{-2}\det A_{11}$, and $\mathbf t\in\mathcal R$ if and only if $\det G(\mathbf t)\ne 0$.

Now fix $b$ and put $\Gamma=\Gamma_{t_b}$. For $\omega=(b,x)\in\{b\}\times\Gamma$ and any other row or column index $\nu$ of $G(\mathbf t)$, the entries of $G(\mathbf t)$ at the positions $(\omega,\nu)$ and $(\nu,\omega)$ do not depend on $x$: they are $\eta(b,b')$ and $\eta(b',b)$ if $\nu\in\{b'\}\times\Gamma_{t_{b'}}$, and they are $(f_j)_b$ if $\nu$ is the $j$-th of the last $c$ indices. So Lemma~\ref{lem:oneblock} and Corollary~\ref{cor:forbidden} apply with $I_0=\{b\}\times\Gamma$, $u=\one$ and $A=C_\Gamma(\psi+\eta(b,b))$, and the matrix $R$ and the vectors $y,z$ there depend only on $b$ and on the entries $t_{b'}$ with $b'\ne b$; let $\theta$ be as in Corollary~\ref{cor:forbidden}. As $A+t\one\one^{\mathsf T}=C_\Gamma(\psi+\eta(b,b)+t)$, the pair $(A,\one)$ is live if and only if $t_b$ is $\psi$-live, and then $f(A,\one)=f_\psi(t_b)-\eta(b,b)$, where $\infty-\eta(b,b)=\infty$. Hence the hypotheses of Lemma~\ref{lem:blockrepl} hold for the $\psi$-live types, the forbidden values $f_\psi$, and $\theta_b(\mathbf t)=\theta+\eta(b,b)$, where $\infty+\eta(b,b)=\infty$ and $\ast+\eta(b,b)=\ast$. If $\psi$ is constant, then by Lemma~\ref{lem:circulant} only the $1$-element sets are $\psi$-live; as every non-constant tuple has an entry that is not $\psi$-live, $\mathcal R=\emptyset$, a contradiction. Otherwise, there are at least three $\psi$-live types by Lemma~\ref{lem:circulant}, and Lemma~\ref{lem:blockrepl} yields a contradiction.

\emph{Case 2b: $c=p$.} Then $N=K$ by \eqref{eq:Jc}, so $W_{\mathbf t}=\bigoplus_bW_{\mathbf t,b}$ with $W_{\mathbf t,b}=K_b\cap\FF_p(\{b\}\times\Gamma_{t_b})$. For each type $\tau$ fix $x_\tau\in\Gamma_\tau$, and use the basis of $W_{\mathbf t,b}$ consisting of the vectors $e_{(b,x)}-e_{(b,x_{t_b})}$ with $x\in\Gamma_{t_b}\setminus\{x_{t_b}\}$. Let $A(\mathbf t)$ be the matrix of the values of $\Phi$ on the resulting basis of $W_{\mathbf t}$, ordered block by block; so $\mathbf t\in\mathcal R$ if and only if $\det A(\mathbf t)\ne 0$. By invariance under $\hat h$, there is a bilinear map $\kappa$ on $J^1(\Omega_1)$ such that $\Phi(u,w)=\kappa(\bar u,\bar w)$ for all $b$ and all $u,w\in K_b$, where $\bar u=\sum_xu_{(b,x)}e_x$. For a type $\tau$ with $|\Gamma_\tau|\ge 2$, let $K_\tau$ be the matrix of the values of $\kappa$ on the vectors $e_x-e_{x_\tau}$ with $x\in\Gamma_\tau\setminus\{x_\tau\}$, let $l_\tau$ be the column vector $(x-x_\tau)_{x\in\Gamma_\tau\setminus\{x_\tau\}}$, all of whose entries are nonzero, and call $\tau$ live if $(K_\tau,l_\tau)$ is live, with forbidden value $f(\tau)=f(K_\tau,l_\tau)$. We call the type of the $1$-element sets live with forbidden value $\infty$.

By Lemma~\ref{lem:forms}(b), the block of $A(\mathbf t)$ with rows indexed by the basis of $W_{\mathbf t,b}$ and columns indexed by the basis of $W_{\mathbf t,b'}$, for $b\ne b'$, is $\gamma(b,b')\,l_{t_b}l_{t_{b'}}^{\mathsf T}$ for some $\gamma(b,b')\in\FF_p$. Fix $b$, let $R$ be the submatrix of $A(\mathbf t)$ with rows and columns outside the block $b$, let $y$ and $z$ be the vectors with the entries $\gamma(b,b')(l_{t_{b'}})_x$ and $\gamma(b',b)(l_{t_{b'}})_x$, indexed like the columns of $R$, and define $d$ and $e$ as in Lemma~\ref{lem:oneblock}(b) and $\theta_b(\mathbf t)$ as in Corollary~\ref{cor:forbidden}; these data depend only on $b$ and on the entries $t_{b'}$ with $b'\ne b$. If $|\Gamma_{t_b}|\ge 2$, then Lemma~\ref{lem:oneblock} and Corollary~\ref{cor:forbidden} apply with $u=l_{t_b}$ and $A=K_{t_b}$. If $|\Gamma_{t_b}|=1$, then $W_{\mathbf t,b}=0$ and $\det A(\mathbf t)=d$, which is nonzero if and only if $\theta_b(\mathbf t)\notin\{\infty,\ast\}$. In both cases, the hypotheses of Lemma~\ref{lem:blockrepl} are satisfied. If $\Gamma_\tau=\{x_\tau,x_\tau+k\}$ with $k\ne 0$, then $K_\tau$ is a $1\times 1$ matrix and $l_\tau=(k)$, so $\det(K_\tau+t\,l_\tau l_\tau^{\mathsf T})=K_\tau+k^2t$ with $k^2\ne 0$, and $\tau$ is live. Hence the type of the $1$-element sets and the $(p-1)/2\ge 2$ types of $2$-element sets are live, and Lemma~\ref{lem:blockrepl} yields a contradiction.
\end{proof}

\begin{remark}\label{rem:pge5}
\begin{enumerate}[(a)]
\item Apart from the description of the largest orbits of $P_4$ and of their centralizers in Section~\ref{sec:orbits}, the proof of Theorem~\ref{thmE} only uses that $p$ is a prime with $p\ge 5$. For $p=3$, Propositions~\ref{prop:Mforced}, \ref{prop:criterion} and~\ref{prop:matrix} remain valid; in particular, Proposition~\ref{prop:Mforced} explains why the computation in the proof of Theorem~\ref{thm:m4} found exactly one candidate for $T\cap H_3$ in each case. But Lemma~\ref{lem:circulant} fails for $p=3$, since a block of three points has only two types, and Lemma~\ref{lem:blockrepl} needs three live types. This is in accordance with Corollary~\ref{cor:TUTD}(c) and Theorem~\ref{thm:m4}.
\item Lemma~\ref{lem:forms}(a) does not extend to $c=p$: the restriction of a $P_2$-invariant bilinear map on $J^p(\Omega_2)=K$ to $K_b\times K_b$ need not be given by a circulant matrix, as the bilinear map $(u,w)\mapsto\langle\tilde u,w\rangle$, where $\tilde u\in\FF_pB_b$ satisfies $\tilde uh_b-\tilde u=u$, shows. This is why Case~2b of the proof is treated separately.
\end{enumerate}
\end{remark}

\section{Concluding remarks}\label{sec:concl}

For $p=3$, Corollary~\ref{cor:TUTD} shows that each of the two largest orbits of $P_m$ on $V_m$ is the unique regular orbit of a subgroup of $P_m$ of order $3^{x_m}$, for every $m\ge 2$. However, by Proposition~\ref{prop:nonunique} and Theorem~\ref{thm:m3}, such a subgroup is in general not determined by its regular orbit, and the number of conjugacy classes of such subgroups grows rapidly with $m$: it is $2$, $30$ and $4914$ for $m=2,3,4$, respectively (Corollary~\ref{cor:m2} and Theorems~\ref{thm:m3} and~\ref{thm:m4}). A conceptual description of these groups for all $m$, and the determination of the isomorphism types of the pairs $(T,V_m)$ for $m\ge 4$, remain open. It would also be interesting to explain why, for $m=3$ and $m=4$, the two largest orbits are the unique regular orbit of the same number of conjugacy classes of subgroups.

For Mersenne primes $p\ge 7$ the situation is quite different. By Theorems~\ref{thmC} and~\ref{thmE}, subgroups of $P_m$ of order $p^{x_m}$ with exactly one regular orbit exist if and only if $m\le 3$; for $m=3$, the groups $T_\Delta$ with $|\Delta|\in\{1,p-1\}$ are irreducible examples (Remark~\ref{rem:generalp}). We have not classified all such subgroups for $m=3$. Moreover, Theorem~\ref{thmE} does not exclude subgroups of $P_m$ of order less than $p^{x_m}$ with exactly one regular orbit for $m\ge 4$.

At the beginning of this work, the authors hoped to classify all $p$-groups having a unique regular orbit on a finite faithful irreducible module of coprime characteristic. This problem remains open, including the case of regular orbits that are not of the largest possible size and the case $p=2$.
\vspace{0.5cm}
\begin{center}
    \textbf{Acknowledgement}
\end{center}
 The first author was supported by a grant from the Simons Foundation (SFI-MPS-TSM-00014153 to Thomas M. Keller).
 \vspace{0.5cm}
 
\begin{center}
    \textbf{AI statement}
\end{center}

This project was started a few years ago with the authors working on and off the project. The onset of AI now allowed the authors to finiths the paper and get more complete answers to their questions than they otherwise would have been able to obtain. We used Claude by Anthropic to correct and add to our draft. In particular, Claude performed all the computational work and obtained the results in Section 5 which answered our original question. After using this tool, the authors reviewed and edited the content as needed and take full responsibility for the final mathematical and narrative integrity of the paper.

\end{document}